\documentclass[11pt]{amsart}
\usepackage{amsmath,amscd,amsfonts,amssymb,epsf,latexsym, eucal}
\newtheorem{theorem}{Theorem}
\newtheorem{lemma}{Lemma}
\newtheorem{corollary}{Corollary}
\newtheorem{proposition}{Proposition}
\theoremstyle{definition}

\newtheorem{remark}{Remark}
\newtheorem{example}{Example}

\def\E{{\mathcal E}}
\def\F{{\mathcal F}}
\def\I{{\mathcal I}}

\def\O{{\mathcal O}}
\def\L{{\mathcal L}}
\def\M{{\mathcal M}}

\def\det{\mathop{\rm det}}

\def\Ext{\mathop{\rm Ext}}
\def\Hom{\mathop{\rm Hom}}
\def\Pic{\mathop{\rm Pic}}
\def\dim{\mathop{\rm dim}}
\def\rank{\mathop{\rm rank}}

\def\Im{\mathop{\rm Im}}

\newcommand{\C}{{\mathbb C}}
\newcommand{\Z}{{\mathbb Z}}
\newcommand{\tensor}{\otimes}
\newcommand{\bP}{{\mathbb P}}
\newcommand{\cal}{\mathcal}
\renewcommand{\H}{\mathcal H}
\newcommand{\Spec}{{\rm Spec}\,}

\newcommand{\G}{\mathcal G}

\newcommand{\norml}{\left\|}
\newcommand{\normr}{\right\|}

\newcommand{\fr}[2]{\frac{#1}{#2}}
\newcommand{\alp}{\alpha}

\newcommand{\eps}{\epsilon}

\newcommand{\rarr}{\rightarrow}

\newcommand{\R}{{\mathbb R}}

\newcommand{\boun}{\partial}

\newcommand{\transv}{\frown\!\!\!\!\mid\,}

\begin{document}
\title{Diagonal subschemes and vector bundles}
\thanks{Pragacz's research supported by a KBN grant and by the Humboldt
Stiftung during his stay at the MPIM in Bonn in the period May-July 2006.}
\author{Piotr Pragacz}
\address{Institute of Mathematics of Polish Academy of Sciences\\
 \'Sniadeckich 8, 00-956 Warszawa, Poland}
\email{ P.Pragacz@impan.gov.pl}
\author{Vasudevan Srinivas}
\address{ School of Mathematics, Tata Institute of Fundamental Research\\
Homi Bhabha Road, Colaba, Mumbai-400005, India}
\email{srinivas@math.tifr.res.in}
\author{Vishwambhar Pati}
\address{Stat-Math Unit, Indian Statistical Institute\\
RV College Post, 8th Mile, Mysore Road, Bangalore-560059, India}
\email{pati@isibang.ac.in}



\date{15.11.2005; revised 29.12.2006}

\begin{abstract}
We study when a smooth variety $X$, embedded
diagonally in its Cartesian square, is the zero scheme of a section of
a vector bundle of rank $\dim(X)$ on $X\times X$. We call this the {\em
diagonal property} (D). It was known that it holds for all flag manifolds
${\rm SL}_n/P$.

We consider mainly the cases of proper smooth varieties, and the
analogous problems for smooth manifolds (``the
topological case'').

Our main new observation in the case of proper varieties is a relation
between (D) and cohomologically trivial line bundles on $X$, obtained by
a variation of Serre's classic argument relating rank 2 vector bundles
and codimension 2 subschemes, combined with Serre duality. Based on this,
we have several detailed results on surfaces, and some results in higher
dimensions.

For smooth affine varieties, we observe that for an affine algebraic
group over an algebraically closed  field, the diagonal is in fact a
complete intersection; thus (D) holds, using the trivial bundle. We
conjecture the existence of smooth affine complex varieties for which (D)
fails; this leads to an interesting question on projective modules.

The arguments in the topological case have a different flavour, with
arguments from homotopy theory, topological K-theory, index theory etc.
There are 3 variants of the diagonal problem, depending on
the type of vector bundle we want (arbitrary, oriented or complex).
We obtain a homotopy theoretic reformulation of the diagonal property as
an extension problem for a certain homotopy class of maps. We also have
detailed results in several cases: spheres, odd dimensional complex
projective quadric hypersurfaces, and manifolds of even dimension
$\leq 6$ with an almost complex structure.

\end{abstract}

\maketitle

\tableofcontents

\bigskip

\centerline{\it To Professor Jean-Pierre Serre on his 80th birthday}

\bigskip

\section{Introduction}

{\it Diagonal subschemes} are important in many questions
of intersection theory. Apart from classical ``reduction to the diagonal''
(cf., e.g., \cite{S2} or \cite{F2}), it was shown in \cite{P}, Sect.5 that
knowing the fundamental class of the diagonal of a variety, is an 
important step towards computing the fundamental classes of {\it all} 
subschemes of this variety (see also \cite{PR}, \cite{Gr}).

A good resolution of the structure sheaf of the diagonal over the 
structure sheaf of $X\times X$ has been used to give a description of the 
derived category $D(X)$ of $X$ \cite{K}, and has proved useful in studying
algebraic $K$-theory of homogeneous spaces and their twisted forms (cf. 
\cite{LSW}, \cite{B}).

We recall an interesting case, where the diagonal is described in a 
suitable fashion, leading to a resolution of its structure sheaf by a 
Koszul complex. Let $G=G_r(V)$ be the Grassmannian parametrizing all 
$r$-subspaces of a vector space $V$ (so, in particular, for $r=1$, we 
consider projective spaces). The Grassmannian is endowed with the 
``tautological" sequence of vector bundles
$$
0 \to S \to V_G \to Q \to 0\,,
$$
where $\rank(S)=r$. Let $G_1=G_2=G$ and use analogous notation for the
tautological vector bundles on $G_i$, $i=1,2$. Denote by
$$
p_1, p_2: G_1\times G_2 \to G
$$
the two projections. Then the diagonal of $G$ is the zero scheme of the
section $s$ of the bundle
$$
\underline{\Hom}(p_1^*S_1, \ p_2^*Q_2)
$$
of rank equal to $\dim(G)$ on $G_1\times G_2$, where $s$ is induced by
the following vector bundle homomorphism:
$$
p_1^*S_1\to p_1^*V_{G_1}=V_{G_1\times G_2}=p_2^*V_{G_2}\to p_2^*Q_2\,.
$$
(This was surely observed independently by so many people that the
``paternity" is impossible to detect.\footnote{-- even by analyzing the 
DNA.})

\smallskip

Let $X$ be a smooth variety. Denote by $\Delta\subset X \times X$ the 
diagonal subscheme, that is, the image of the diagonal embedding
$$
\delta: X \hookrightarrow X \times X \,,
$$
given by $\delta(x)=(x,x)$.

\smallskip

We are interested in when the following {\it diagonal property} holds:

\medskip

``There exists a vector bundle $\E$ of rank equal to $\dim(X)$ on $X\times 
X$
and a section $s$ of $\E$, such that $\Delta$ is the zero scheme of $s$.''

\smallskip

\noindent
In the following, we shall use abbreviation ``(D)'' for the diagonal 
property. 

Observe that if $X_1$ and $X_2$ satisfy (D), then it also holds 
for $X_1\times X_2$. Moreover, (D) is obviously valid for curves. We 
noted that it holds for any Grassmannian; it also holds for flag varieties 
of the form $SL_n/P$ over any field (cf. \cite{F}, \cite{FP}).
In fact, in this last case, (D) was a starting point for the development
of the Lascoux-Sch\"utzenberger theory of {\it Schubert polynomials} 
for $SL_n$ \cite{LS}
(see also \cite{F}, \cite{FP}). The question: 

\medskip

``Do the flag varieties for other classical groups have (D)?''

\smallskip

\noindent
arose in discussions of the first author with William Fulton while writing up \cite{FP}
at the University of Chicago in 1996.

\bigskip

In the present paper, we investigate this property mainly for surfaces,
and have some results in the higher dimensional case. We use an argument
arising from the fundamental {\em Serre construction} relating
codimension 2 subschemes and vector bundles of rank $2$ (cf. \cite{S1}),
combined with {\em Serre duality}, to relate (D) with the  existence or
absence of {\it cohomologically trivial line bundles}.

In the surface case, this leads to the following result, summarizing
several conclusions.

\begin{theorem}\label{th-surf}
Let $X$ be a smooth projective surface over an algebraically closed field.
\begin{enumerate}
\item[(a)] There exists a birational morphism $f:Y\to X$ such that $Y$
satisfies (D).
\item[(b)] If $Y\to X$ is a birational morphism, $X$ satisfies (D), and
$\Pic(X)$ is finitely generated, then $Y$ satisfies (D).
\item[(c)] Suppose $X$ is birational to one of the following: a ruled or
an abelian surface, or a K3 surface with two disjoint rational curves, or
an elliptic fibration with a section, or a product of 2 curves, or
 a complex Enriques or hyperelliptic surface. Then $X$ satisfies (D).
\item[(d)] Suppose $\Pic(X)=\Z$, such that the ample generator of
$\Pic(X)$ has a nonzero section, and $X$ satisfies (D). Then $X\cong
\bP^2$. In particular, (D) fails for general algebraic K3 surfaces, or for
general hypersurfaces $X\subset\bP^3$ of degree $\geq 4$.
\end{enumerate}
\end{theorem}

For higher dimensional varieties, the general result we have is the
following.
\begin{theorem} Let $X$ be a smooth projective variety with $\Pic(X)=\Z$,
such that the ample generator $\O_X(1)$ has a nonzero section. Suppose $X$
has (D).  Then $\omega_X\cong\O_X(-r)$ for some $r\geq 2$, i.e., $X$ is a
Fano variety of index $\geq 2$.
\end{theorem}

On the other hand, (D) {\it fails} for a smooth projective quadric of any
odd dimension $\geq 3$. Curiously, our proof of this for quadrics of
dimension $\geq 5$ is by reducing to the case of quadrics over the complex
numbers, and using the topological results below.

We have not been able to decide if (D) holds for cubic 3-folds, though we
suspect it does not hold, at least in general.

While studying (D), often, a related {\it point property}
(cf. Section \ref{ptpr}) is particularly useful (both in disproving
and proving (D)). The point property is closely related to the property
for a point to be a {\it complete intersection}; this property was
extensively studied in \cite{Mur}. However, it seems likely that it is a
strictly weaker property: we can verify it for a cubic 3-fold, for
example.

(D) makes sense for smooth affine varieties as well. Here, Serre's
construction implies that (D) holds for all affine surfaces. We conjecture
that there exist smooth affine complex 3-folds for which (D) fails.
We formulate a question on projective modules, a negative answer to which
gives examples of smooth affine varieties for which (D) would fail (see
Section~5).

There are also topological versions of (D). Here, the question comes in
several flavours: we may ask if the property holds
\begin{enumerate}
\item[(i)] for a smooth manifolds, with an arbitrary smooth vector bundle
of appropriate rank (this property is denoted $(D_r)$), or
\item[(ii)] for oriented manifolds, and oriented vector bundles (this is
denoted $(D_o)$), or
\item[(iii)] (in appropriate cases) for an even dimensional almost
complex smooth manifold $X$, and a compatible complex vector bundle on
$X\times X$ (this is called property $(D_c)$).
\end{enumerate}
Some of these  questions are investigated in this paper, in Section~6. 

There are some obvious remarks, analogous to those in the 
algebraic case: each of the diagonal properties is compatible with 
products, $(D_r)$ holds for 1-manifolds, and real flag varieties 
$SL_n(\R)/P$, and $(D_c)$ holds for Riemann surfaces, and for flag 
varieties $SL_n(\C)/P$. However, orientable real flag varieties 
need not satisfy $(D_o)$ (it fails for odd dimensional real projective 
spaces, as we see in Theorem~\ref{top}). 

We obtain a homotopy theoretic translation of the diagonal property (of 
any of the above types): if $Y$ is the complement of an open tubular
neighborhood of the diagonal in $X\times X$, then the property is
equivalent to the extendability of a certain vector bundle on $\partial
Y$ to the whole of $Y$; this can be viewed as an extension problem 
for a map to a classifying space, leading to a possible obstruction 
theoretic approach (see Lemma~\ref{lemma4}).

We also have several explicit results, obtained by different methods, 
which are as follows.
\begin{theorem}\label{top}
\begin{enumerate}
\item The sphere $S^n$ has $(D_r)$ iff $n=1,2,4$ or $8$. It has $(D_o)$ 
iff $n=2,4,8$. It has $(D_c)$ iff $n=2$.
\item A compact oriented odd dimensional manifold $X$ does not have 
$(D_o)$, and does not have $(D_r)$ if $H^1(X,\Z/2\Z)=0$.
\item Any compact almost complex $4$-manifold has $(D_c)$. If a compact 
almost complex $6$-manifold $X$ satisfying
$H^1(X,\Z)=0$ and $H^2(X,\Z)\cong \Z$ has $(D_c)$, then it has a spin structure.
\item A (smooth) complex quadric projective hypersurface $X\subset
\bP^{2n}(\C)$ does not have $(D_c)$, unless $n=1$.
\end{enumerate}
\end{theorem}

There are other contexts in which analogues of (D) make sense, that we do 
not investigate here: e.g. for compact complex manifolds, and for 
Stein manifolds, or over non-algebraically closed fields. Thus, one 
may ask if there exist Stein manifolds of dimension $\geq 3$ for which 
the diagonal is not the zero set, with multiplicity 1, of a section of a 
holomorphic vector bundle. Similarly, for a projective smooth surface over 
a number field, which is geometrically rational, the validity of (D) 
might provide a subtle obstruction to rationality over the given field. 

\section{Preliminaries on (D)}
\subsection{Some general remarks}

Let $\I_{\Delta}\subset \O_{X\times X}$ denote the ideal sheaf of
$\Delta\subset X\times X$.
If a variety $X$ admits (D) then the cotangent sheaf is isomorphic
to the restriction of $\E^*$ to the diagonal:
\begin{equation}\label{con}
\E^*_{|\Delta}=\I_{\Delta}/\I_{\Delta}^2\cong \Omega_X^1
\end{equation}
(via the isomorphism $\Delta \cong X$),
so that it is locally free and consequently $X$ is {\it smooth}.

Let $\L=\det(\E^*)$, then
\begin{equation}\label{omega}
\L_{|\Delta}=\omega_{\Delta}\cong \omega_X\,,
\end{equation}
We have also the following expression for the fundamental class
of $\Delta$:
\begin{equation}
[\Delta]=c_{\dim(X)}(\E)
\end{equation}
by the Grothendieck formula \cite{G}, Th\'eor\`eme 2.

\subsection{(D) for proper smooth varieties}
In this section, we assume that $X$ is a proper smooth variety of 
dimension $n\geq 2$ over an algebraically closed field $k$. We make a 
preliminary analysis of the condition (D), motivated by Serre's 
construction of vector bundles of rank $2$ (cf. \cite{S1}), and some 
variants.

\smallskip

Consider the exact sequence
\begin{equation}\label{ex}
0 \to \I_{\Delta} \to \O_{X\times X} \to \O_{\Delta} \to 0\,.
\end{equation}
For any line bundle $\L$ on $X\times X$, if we apply the functor 
$\Hom(-,\L)$ (and its derived functors) to the exact sequence 
(\ref{ex}), we get the following exact sequence of global $\Ext$'s
(all $\Ext$'s without subscripts are taken over $X\times X$):
\begin{equation}\label{ls2}
H^{n-1}(X\times X, \L) \to {\Ext}^{n-1}(\I_{\Delta}, \L) 
\stackrel{\beta}{\to}
{\Ext}^n(\O_{\Delta}, \L) \stackrel{\alpha}{\to}
H^n(X\times X, \L)\,.
\end{equation}

Suppose now that $\E$ is a vector bundle (i.e., a locally free sheaf) of 
rank $n$ on $X\times X$, with a global section $s:\O_{X\times X}\to \E$, 
whose zero scheme is the diagonal $\Delta$.  We then have the 
(truncated) Koszul complex associated with the dual $s^*:\E^*\to 
\O_{X\times X}$,
\begin{equation}\label{Kosz}
0 \to \det(\E^*) \to \cdots\to \E^* \stackrel{s^*}{\to}
\I_{\Delta} \to 0 \,,
\end{equation}
which determines an element 
$[s]\in \Ext^{n-1}(\I_{\Delta},\det(\E^*))$. On splicing this with the 
short exact sequence (\ref{ex}), we obtain the Koszul resolution of 
$\O_{\Delta}$, and a resulting element
\[\beta([s])\in {\Ext}^n(\O_{\Delta},\det(\E^*)),\]
where $\beta$ is the map in (\ref{ls2}), with the choice $\L=\det(\E^*)$.
This amounts to the assertion that the boundary map $\beta$ may be viewed 
as a Yoneda product with the class of the extension (\ref{ex}).  

Our basic criterion for testing if (D) holds is based on the following 
result, which is folklore (see the discussion in \cite{EPW}, where this is 
the condition that the subscheme $\Delta\subset X\times X$ is ``strongly 
subcanonical'').
\begin{proposition}\label{prop} Let $\L$ be a line bundle on $X\times 
X$ whose restriction to $\Delta$ coincides with $\omega_{\Delta}$. Assume 
that there exist a rank $n$ vector bundle $\E$ on $X\times X$ with 
$\det(\E^*)=\L$, and a section $s\in \Gamma(X\times X,\E)$ satisfying 
$\Im(s^*)=\I_{\Delta} \,$.  Then the map $\alpha$ in the exact sequence 
(\ref{ls2}) vanishes. The converse holds if $n=\dim(X)=2$.
\end{proposition}
\begin{proof}
By the assumption, $\L$, restricted to $\Delta$, coincides with 
$\omega_{\Delta}$. By Serre duality on $X\times X$,
\begin{equation}\label{ser}
{\Ext}^n(\O_{\Delta},\L)^*\cong H^n(\Delta,
\L^* \otimes {\omega_{X\times X}}_{|\Delta})
=H^n(\Delta, \omega_{\Delta})\,.
\end{equation}
Hence the space  ${\Ext}^n(\O_{\Delta},\L)$ is 1-dimensional.

We also have isomorphisms of sheaf Ext's
\begin{equation}
\underline{\Ext}^i(\O_{\Delta}, \L)
\cong\left\{\begin{array}{cl}0&\mbox{if $i\neq n$}\\
\underline{\Hom}_{\O_\Delta}(\det(\I_\Delta/\I_{\Delta}^2),
\L_{|\Delta})&\mbox{if $i=n$}\end{array}
\right.
\end{equation}
(cf. \cite{Ha}, III, \S7) because $\Delta\subset X\times X$ is a
nonsingular subvariety of codimension $n$, with
\begin{equation}
\underline{\Hom}_{\O_\Delta}(\det(\I_\Delta/\I_{\Delta}^2),
\L_{|\Delta})\cong
\omega_{\Delta}^{-1} \otimes \L_{|\Delta} \cong \O_{\Delta}\,.
\end{equation}
This implies that the canonical map
\begin{equation}
\gamma: {\Ext}^{n}(\O_{\Delta}, \L) \to
H^0(X\times X,{\underline{\Ext}}^{n}(\O_{\Delta},
\L))\cong H^0(\Delta,\O_{\Delta})
\end{equation}
is an isomorphism of 1-dimensional vector spaces.

Now, if there exists a vector bundle $\E$ with $\det(\E^*)=\L$,
and a section $s\in \Gamma(X\times X,\E)$ with zero scheme equal to the
diagonal, then there exists $\beta[s]\in {\Ext}^{n}(\O_{\Delta}, \L)$
which is in the kernel of
$\alpha$. But $\beta[s]$ is the class of the Koszul resolution of
$\O_{\Delta}$, which is a generator of the 1-dimensional vector space
$H^0(X\times X,{\underline{\Ext}}^{n}(\O_{\Delta}, \L))$. Hence $\alpha$,
which has a non-trivial kernel, and a 1-dimensional domain,  must be the
zero map.

When $n=2$, if $\alpha=0$, the isomorphism $\gamma$ implies the existence
of an element $[s]\in \Ext^1(\I_{\Delta},\L)$ whose image $\beta([s])\in
\Ext^2(\O_{\Delta},\L)$ generates this 1-dimensional vector space.
This element $[s]$ determines an extension,  that is to say, a vector
bundle $\E$ and section $s$ giving rise to a 3-term exact sequence
(\ref{Kosz}); this is the Serre construction.
\end{proof}

We call a line bundle $\M$ on $X$ {\em cohomologically trivial} if
$H^i(X,\M)=0$ for all $i$. If $\M$ is cohomologically trivial on a smooth
proper variety $X$, then by Serre duality, so is $\M^{-1}\tensor\omega_X$.

We now discuss the criterion which we actually use in studying the
property (D), in many instances. This leads to a good understanding of (D)
in the surface case, and a nontrivial necessary condition for (D) to hold,
in the higher dimensional case.

\begin{theorem}\label{Tm}
(i) Denote by $p_1, p_2: X \times X \to X$ the two projections.
Suppose that
\begin{equation}\label{prod}
\Pic(X\times X)\cong p_1^*\Pic(X)\oplus p_2^*\Pic(X)\,,
\end{equation}
and that $X\times X$ supports a vector bundle $\E$ with a section $s$
such that (D) holds. Then there exists a cohomologically trivial line
bundle $\M$ on $X$, such that $\det(\E)=p_1^*\M^{-1}\tensor
p_2^*(\M\tensor\omega_X^{-1})$.

\noindent
(ii) If $\dim(X)=2$, and there exists a cohomologically trivial line bundle
on $X$, then (D) holds for $X$.
\end{theorem}
\begin{proof}
Suppose that there exists a vector bundle $\E$ on $X\times X$
of rank $n$ such that the diagonal is the zero scheme of its section $s$.
Let $\L=\det(\E^*)$, and form the corresponding exact sequence
(\ref{ls2}). From Proposition~\ref{prop}, we must have $\alpha=0$.

Now consider the dual linear transformation to $\alpha$:
\begin{equation}\label{pair}
\alpha^*: H^n(X\times X,\L)^* \to {\Ext}^n(\O_{\Delta}, \L)^*
\cong H^n(\Delta, \omega_{\Delta})=k\,.
\end{equation}
Using (\ref{prod}) and (\ref{omega}), choose $\M \in \Pic(X)$ such that
$$
\L=\det(\E^*)\cong p_1^*(\M) \otimes p_2^*(\M^{-1}\otimes \omega_X)\,.
$$
By Serre duality on $X\times X$, we get that
\[H^n(X\times X,\L)^*\cong H^n(X\times X,\L^{-1}\tensor\omega_{X\times
X})\cong H^n(X\times X,p_1^*(\M^{-1}\tensor\omega_X)\tensor p_2^*\M).\]
From the K\"unneth formula, we have
\begin{equation}\label{kunneth}
H^n(X\times X, p_1^*(\M^{-1}\tensor\omega_X) \otimes p_2^*(\M))
=\oplus_{i=0}^n H^i(X,\M^{-1}\tensor\omega_X) \otimes H^{n-i}(X,
\M)\,.
\end{equation}
Further, on any summand on the right, the induced map
$$H^i(X,\M^{-1}\tensor\omega_X) \otimes H^{n-i}(X, \M)\hookrightarrow
H^n(X\times X,
p_1^*(\M^{-1}\tensor\omega_X) \otimes p_2^*(\M))\stackrel{\alpha^*}{\to}
H^n(\Delta,\omega_{\Delta})=k$$
coincides with the {\em Serre duality pairing} on cohomology of $X$, and
is hence a non-degenerate bilinear form, for each $0\leq i\leq n$.

Thus, the dual map to $\alpha$ vanishes if and only if all of the summands
on the right side of (\ref{kunneth}) vanish, which amounts to saying that
$\M$ is cohomologically trivial.

Conversely, if $\M$ is cohomologically trivial, then in the exact
sequence (\ref{ls2})
determined by the line bundle
\[\L=p_1^*\M\tensor  p_2^*(\M^{-1}\tensor\omega_X),\]
the map $\alpha$ is the zero map, by reversing the above argument. Hence,
if $n=2$, we deduce that $X\times X$ supports a vector bundle $\E$ of
rank $2$ and a section $s$ with zero scheme $\Delta$, by the surface case
of Proposition~\ref{prop}.
\end{proof}

In particular, if the isomorphism (\ref{prod}) holds for $X$, and $X$ supports
{\it no} cohomologically trivial bundle, then $X$ fails to have (D).

\smallskip

In order to apply Theorem~\ref{Tm} in various situations, we need to be
able to verify the hypothesis (\ref{prod}) on the Picard group. This
property is well understood, and the facts are recapitulated below.

\begin{lemma}\label{Pic} Let $X$ be any smooth proper variety over an
algebraically closed field $k$.
\begin{enumerate}
\item[(i)]
The isomorphism (\ref{prod}) holds $\Leftrightarrow$ $\Pic(X)$
is a finitely generated abelian group.
\item[(ii)] If \ $H^1(X, {\cal O}_X)=0$, then (\ref{prod}) holds.
\end{enumerate}
\end{lemma}
\begin{proof} (i) \ Let $\underline{\rm Pic}^0(X)$ denote the neutral component
(identity component) of the {\em Picard scheme} of $X$, in the sense of
Grothendieck. Then $\underline{\rm Pic}^0(X)$ is a projective connected
$k$-group scheme, such that
\begin{enumerate}
\item[(a)] its Lie algebra is naturally identified with the vector space
$H^1(X, {\mathcal O}_X)$\,;
\item[(b)]  the associated reduced group scheme is an abelian variety,
the classical {\em Picard variety} (in the sense of Weil).
\end{enumerate}

The group of connected components of the Picard scheme is identified with
the N\'eron-Severi group, which is always a finitely generated abelian
group. The group of $k$-rational points of a positive dimensional abelian
variety is not finitely generated. Thus, $\Pic(X)$ is a finitely
generated group precisely when $\underline{\rm Pic}^0(X)$ is
a 0-dimensional, local $k$-group scheme (in characteristic
0, this means it is 0). The connected group scheme  $\underline{\rm
Pic}^0(X)$ is of course 0 if its Lie algebra $H^1(X,{\mathcal O}_X)$
vanishes.

In any case, if $\L$ is any line bundle on $X\times X$, which (for some
fixed base point $x_0\in X$) is trivialized on $X\times\{x_0\}$ as well
as on $\{x_0\}\times X$, then by the universal property of the Picard
scheme, $\L$ determines a morphism of $k$-schemes $f_{\L}:X\to
\underline{\rm Pic}(X)$, which maps $x_0$ to the identity point, so that
$\L$ is the pull-back under $1_X\times f_{\L}:X\times X\to X\times
\underline{\rm Pic}(X)$ of a suitable Poincar\'e bundle. Clearly $f_{\L}$
factors through the neutral component, and then through the
corresponding reduced scheme, which is the identity point, when $\Pic(X)$
is finitely generated. Hence $\L$ is the trivial bundle. This clearly
implies (\ref{prod}).

On the other hand, suppose $\Pic(X)$ is not finitely generated, which is
to say that $\underline{\rm Pic}^0(X)_{\rm red}$ is a positive
dimensional abelian variety. We claim that (\ref{prod}) {\em does
not} hold in this case.

Indeed, the dual abelian variety to $\underline{\rm Pic}^0(X)_{\rm red}$
is the Albanese variety ${\rm Alb}(X)$, in the sense of Weil, which is an
abelian variety of the same dimension ($>0$) as the Picard variety;
further, there is a morphism $f:X\to {\rm Alb}(X)$ such that $f(x_0)=0$,
which is universal among morphisms of pairs $(X,x_0)\to (A,0)$ with $A$
any abelian variety. As a consequence, the $k$-points of the image
variety $f(X)$ generate ${\rm Alb}(X)$ as a group, and in particular, $f$
has positive dimensional image.

From the duality theory of abelian varieties, we know that the Albanese
and Picard varietes of $X$, being mutually dual, are also {\em isogenous}.
This means that, composing $f$ with an isogeny, we can construct a
morphism $g:X\to \underline{\rm Pic}^0(X)_{\rm red}$ with
$g(x_0)$ equal to the neutral element, whose image is a positive
dimensional subvariety, whose $k$-points generate the group of
points of the abelian variety $\underline{\rm Pic}^0(X)_{\rm red}$. The
pullback of
the Poincar\'e line bundle under the morphism
\[1_X\times g: X\times X \to X\times \underline{\rm Pic}^0(X)\]
is a line bundle on $X\times X$ which is non-trivial, but is trivial when
restricted to $X\times\{x_0\}$ or to $\{x_0\}\times X$. Hence (\ref{prod}) fails.

\smallskip

\noindent
(ii) This follows from ``Cohomology and Base Change'' (cf., e.g.,
\cite[Theorem 12.11]{Ha} and comment on p. 292).

\end{proof}

\section{Main results for surfaces}

We now explore (D) for surfaces. Our main references on the theory of
surfaces are: \cite{BPV}, \cite{Bea}, \cite{GH}, and \cite{Ha}.
From Theorem~\ref{Tm}, for a surface $X$, (D) is more or less equivalent
to the existence of a cohomologically trivial line bundle on $X$.

Before a more systematic discussion, we consider the following simple
result, since it was our first example of failure of (D).

\begin{proposition}\label{Pk3} Let $X$ be a generic complex algebraic K3
surface. Then (D) fails for $X$.
\end{proposition}
\begin{proof}
It follows from e.g. \cite{GH}, p. 594, \ that $\Pic(X)=\Z\cdot \O_X(1)$,
where $\O_X(1)$ is an ample line bundle on $X$. Let $\O_X(1)\cong \O_X(D)$
for a divisor $D$; then $d:=D^2>0$.

Since $X$ has a trivial canonical bundle and $\chi(X,\O_X)=2$,
the Riemann-Roch theorem implies that for any integer $n$,
\[\chi(X,\O_X(n))=\frac{n^2d}{2}+2\geq 2.\]
Consequently, since every line bundle $\M$ is isomorphic to some $\O_X(n)$,
$\M$ must has some nontrivial cohomology. The assertion now follows
from Theorem \ref{Tm}(i) (and Lemma \ref{Pic}(i)).
\end{proof}

We now proceed more systematically.

\subsection {(D) within a birational class}

A classical result in surface theory asserts that any smooth proper
(equivalently, projective) algebraic surface is obtained from a
{\em (relative) minimal  model} by a succession of blow ups of points
(cf., e.g., \cite{Ha} III, Theorem 5.8 and II, 4.10.2).  Thus, the
following result reduces (D) for surfaces to the question of finding
cohomologically trivial line bundles on relative minimal models.

\begin{proposition}\label{blowup}
Given a birational morphism $f: X\to Y$ of smooth surfaces, if $Y$ supports
a cohomologically trivial line bundle $\M$, then
the pullback $f^*\M$ is also cohomologically trivial.
\end{proposition}
\begin{proof}
The morphism $f$ is a composition of point blow ups, as noted above. By
induction on the number of blow ups, we thus reduce to the case when $f$
is the blow up of 1 point. In this case, it is standard that
$f_*\O_X=\O_Y$ and $R^if_*\O_X=0$ for all $i>0$ (see \cite{Ha}, V, 3.4, for
example). From the projection formula (\cite{Ha}, III, Ex. 8.3, for
example), it follows that for any cohomologically trivial line bundle $\M$
on $Y$, the pullback $f^*\M$ is cohomologically trivial on $X$.
\end{proof}
\begin{corollary} If $X$ admits a morphism to a relative minimal model $Y$
which supports a cohomologically trivial line bundle, then $X$ has (D).
\end{corollary}
Let us note also:
\begin{corollary}
If $f:X \to Y$ is a birational morphism of smooth surfaces with finitely
generated Picard groups, and
$Y$ has (D), then $X$ has (D).
\end{corollary}

On the other hand, we claim that any ``sufficiently non-minimal'' surface
has a cohomologically trivial line bundle, and hence satisfies (D).

\begin{theorem} If $X$ is any smooth projective surface over an
algebraically closed field $k$, there is a birational morphism $f: Y\to X$
of smooth projective surfaces such that $Y$ satisfies (D).
\end{theorem}
\begin{proof}
Choose a line bundle $\L$ on $X$ so that $H^i(X,\L)=0$ for $i>0$ (by
Serre vanishing, $H^i(X,\O_X(n))=0$ for all $i>0$, and any $n>>0$, so
such an $\L$ certainly exists). If $H^0(X,\L)=0$, then $\L$ is
cohomologically trivial, hence (D) holds for $X$ itself. So let us assume
that $\dim H^0(X,\L)=r>0$.

For any $x\in X$, let $\L(x)=\L_x\otimes_{\O_{X,x}} k(x)$ denote the fiber
of $\L$ at $x$; this is a 1-dimensional vector space over $k=k(x)$.

\smallskip

We now state a standard lemma, whose proof is left to the reader.
\begin{lemma}
Let $X$ be a projective variety over an algebraically closed field, and
$\L\in \Pic(X)$ with  $\dim H^0(X,\L)=r>0$. Then for any $r$ general points
$x_1,\ldots ,x_r$ of $X$, the induced map of $r$-dimensional vector spaces
\begin{equation}
f_r: H^0(X,\L) \to \oplus_{i=1}^r \L(x_i)
\end{equation}
is an isomorphism. \footnote{This says that, if we blow up the base scheme
to get a morphism to projective space ${\Bbb P}^{r-1}$, the linear span
of the image is the whole projective space, so that there exist $r$ linearly
independent points in the image.}
\end{lemma}

Now let $\L$ be as before on our surface $X$, and let $x_1,\ldots ,x_r$ be
chosen as in the lemma. Let $f:Y\to X$ be the blow up of the points
$x_1,\ldots, x_r$, and let $E_1,\ldots, E_r$ be the corresponding exceptional
curves. Consider the line bundle
$$
\M=f^*\L\otimes \O_Y(-E_1-E_2-\cdots -E_r)
$$
on $Y$. We claim that $\M\in\Pic(Y)$ is cohomologically trivial, and hence 
$Y$ satisfies (D).
Indeed, it is easy to see that if $S=\{x_1,...,x_r\}$, and $\I_S$
is the ideal sheaf of $S$, then
$$
R^if_*\M = \L\otimes R^if_*\O_Y(-E_1-\cdots -E_r)=0
$$
if $i>0$, and $f_*\M = \L\otimes \I_S$.
Thus, from the Leray spectral sequence for $f$, it follows at once that
$$
H^i(Y,\M)=H^i(X,\L\otimes \I_S)\,.
$$
There is an exact sequence of sheaves
$$
0\to \L\otimes \I_S \to \L \to \L\otimes \O_S \to 0\,,
$$
where we may identify
$$
\L\otimes \O_S= \oplus_{i=1}^r \L(x_i)_{x_i}
$$
(here, for a point $x$ and abelian group $A$, we let $A_x$ denote the
skyscraper sheaf with stalk $A$ at $x$).  Then clearly
$H^i(X,\L\otimes \O_S)=0$ for $i>0$, while
$$
H^0(X,\L)\to H^0(X,\L\otimes \O_S)=\oplus_{i=1}^r \L(x_i)
$$
is the map $f_r$ considered above, which (by the choice of the set
$S$) is an isomorphism.

Hence the long exact cohomology sequence for the above sequence of sheaves
implies that $\L\otimes \I_S$ is cohomologically trivial on $X$,
and so $\M$ is cohomologically trivial on $Y$.
\end{proof}

\subsection{Ruled surfaces}
We shall examine now the case of birationally ruled surfaces.

It is useful for us to make explicit the following observation (which is
of course well-known).
\begin{lemma}\label{ctl-curves}
Any smooth projective curve $Y$ over an algebraically closed field 
supports a cohomologically trivial line bundle $\L$. 
\end{lemma}
\begin{proof}
This is clear if $Y=\bP^1$ (take $\L=\O_Y(-1)$). If $Y$ has genus $g>0$, 
the isomorphism classes of line bundles of degree $g-1$ on $Y$ are 
parametrized by a $g$-dimensional variety $J^{g-1}(Y)$ (it is a 
principal homogeneous space under the Jacobian variety $J(Y)$). 
On the other hand, the subvariety parametrizing line bundles with a 
non-zero section is the image of the natural morphism $S^{g-1}(Y)\to 
J^{g-1}(Y)$, with domain the $(g-1)$st symmetric power of
$Y$, which is the parameter variety for effective divisors of degree 
$g-1$. Clearly this is not surjective; a point in the complement of the 
image corresponds to a line bundle $\L$  of degree $g-1$ on $Y$ such that 
$\L$ has no non-zero global sections. Now the Riemann-Roch theorem on $Y$
implies $\L$ has vanishing Euler characteristic, hence is cohomologically
trivial.
\end{proof}

\begin{proposition}\label{ruled} Let $X$ be a birationally ruled surface.
Then $X$ admits a cohomologically trivial line bundle, and in particular,
(D) holds for $X$.
\end{proposition}
\begin{proof} We may as well assume $X\not\cong\bP^2$, since
$\O_{\bP^2}(-1)$ is cohomologically trivial. Then, from the
classification of birationally ruled surfaces, we know that
$X$ is birational to  $Y\times \bP^1$ for some smooth projective
curve $Y$, such that there is a morphism $\pi:X\to Y$ with general fiber
$\bP^1$. The morphism $\pi$ has a factorization as a composition $X\to
\bar{X}\to Y$, where $X\to \bar{X}$ is a composition of point blow ups,
and  $\bar{X}\to Y$ is a $\bP^1$-bundle. By Proposition~\ref{blowup}, it
suffices to show $\bar{X}$ supports a cohomologically trivial line bundle.

Hence we may without loss of generality assume that $\pi:X\to Y$ is a
$\bP^1$-bundle. We now note that if $\L$ is cohomologically trivial on
$Y$ (and such a line bundle exists, by Lemma~\ref{ctl-curves}), then
$\pi^*\L$  is cohomologically trivial on the surface $X$. This  follows
from the Leray spectral sequence, because for the $\bP^1$-bundle
$\pi:X\to Y$, we have $R^i\pi_*{\cal O}_X=0$ for $i>0$, and
$\pi_*\O_X=\O_Y$.
\end{proof}

\subsection{Surfaces of Kodaira dimension 0}

\begin{proposition} An abelian surface supports a cohomologically
trivial line bundle. Thus (D) holds for abelian surfaces.
\end{proposition}
\begin{proof}
It is well known (cf., e.g., \cite{M}, Sect. 8) that for an abelian
variety $X$, any line bundle $\M$ on $X$ which is non-trivial, but is
algebraically equivalent to 0 (i.e., having a nontrivial class in
$\Pic^0(X)$), is cohomologically trivial.

Note that in the surface case, after showing that $H^0(X,\M)=0$
{\it (loc.cit.)}, one can argue as follows. Since $\omega_X$ is trivial,
we get by Serre duality that $H^2(X,\M^{-1})=0$. But, exchanging the roles
of $\M$ and $\M^{-1}$, this implies that $H^2(X,\M)=0$.
On the other hand, we have by the Riemann-Roch theorem that $\chi(\M)=0$,
so that we also have $H^1(X,\M)=0$.
Applying Theorem \ref{Tm}(ii), the assertion follows.
\end{proof}

\begin{proposition}\label{two} A K3 surface $X$ with two disjoint
smooth rational curves supports a cohomologically trivial line bundle,
and (D) holds for it.
\end{proposition}
\begin{proof} Let $D_1$ and $D_2$ be two disjoint smooth rational
curves on $X$.
Let $\M=\O(D_1-D_2)$. We claim that $\M$ is cohomologically trivial.
By the adjunction formula, the curves $D_1, D_2$ are $-2$ curves. 
We thus have $H^0(X,\M)=0$. By Serre duality, invoking that $\omega_X$ 
is trivial, we get that $H^2(X,\M^{-1})=0$. But, exchanging the roles of $D_1$ 
and $D_2$, this gives $H^2(X,\M)=0$. The Riemann-Roch theorem for a K3 surface
reads (with $\L=\O(D)$):
\begin{equation}
\chi(\L)=\frac{1}{2} D^2 +2\,.
\end{equation}
Hence we get the vanishing $\chi(\M)=0$, and consequently we have $H^1(X,\M)=0$. 
Applying Theorem \ref{Tm}(ii), the assertion follows.
\end{proof}

\begin{corollary} For the Kummer surface, (D) holds.
\end{corollary}

In the rest of this section, we assume the ground field is $\C$, for
simplicity.

We use the term {\em hyperelliptic surface} to mean a complex surface $X$
which is a quotient of a product $E\times C$ by the diagonal action of a
finite abelian group $G$ which acts on each of the factors, where $E$ is
an  elliptic curve on which $G$ acts faithfully by translations, while $C$
is a smooth projective curve of genus $>0$ such that $C/G\cong \bP^1$.
The surface $X=(E\times C)/G$ admits a morphism to $C/G\cong \bP^1$ with
fibers which are elliptic curves. Such a surface has Kodaira dimension
$0$ if $C$ is also elliptic, else it has Kodaira dimension 1.

\begin{proposition}\label{bielliptic}  A hyperelliptic surface (as above)
admits a cohomologically trivial line bundle, and (D) holds for it.
\end{proposition}
\begin{proof}
First observe that since $q:E\times C\to X$ is a finite covering of
complex surfaces, the cohomology of any line bundle on $X$ injects into
the cohomology of its pull-back to $E\times C$. Hence, if we show that
there is some line bundle $\M$ on $X$ whose pull-back to $E\times C$ is
cohomologically trivial, then $\M$ itself is cohomologically trivial.

Since $G$ is a finite group of translations on $E$, the map
$\pi:E\to E/G$ is an isogeny of elliptic curves. We can find a nontrivial
line bundle $\L$ on $E/G$ of degree 0, whose pull-back $\pi$ is a
non-trivial line bundle on $E$, also of degree 0 (i.e., $\L$ considered as
a point of $E/G\cong \Pic^0(E/G)$ is not in the kernel of the dual isogeny
to $\pi$).

In particular, $\pi^*\L$ is cohomologically trivial on $E$.
If $p_1:E\times C\to E$ is the projection, the pullback $p_1^*\pi^*\L$ on
$E\times C$ is again cohomologically trivial, from the K\"unneth formula.

There is a commutative diagram, whose horizontal arrows are quotients mod
$G$, with induced vertical arrow $f$
$$
\begin{array}{ccc}
E\times C &\stackrel{q}{\to} &X\\
p_1\downarrow\quad&&\quad\downarrow f\\
E&\stackrel{\pi}{\to} &E/G
\end{array}
$$
Thus $\M=f^*\L$ has the property that $q^*\M=p_1^*\pi^*\L$ is
cohomologically trivial; hence $\M$ is cohomologically trivial on $X$.
\end{proof}

We shall now discuss the case of Enriques surfaces (surfaces with
$p_a=p_g=0$ whose canonical line bundle has order 2).
\begin{proposition}
Any complex Enriques surface supports a cohomologically trivial line
bundle, and hence has (D).
\end{proposition}
\begin{proof}
First, suppose we have a smooth $-2$ curve $E$ on the Enriques
surface $X$ (i.e., $E\cong \bP^1$, and $E^2=-2$). We then claim
$\O_X(-E)$ is cohomologically trivial (and so the diagonal property
holds). Indeed, we have $\chi(X,\O_X(-E))=0$ by the Riemann-Roch theorem. 
Since $\O_X(-E)$ is the ideal sheaf of a smooth rational curve, and $X$ 
has $H^1(X,\O_X)=0$, we have
\begin{equation}
H^0(X,\O_X(-E))=H^1(X,\O_X(-E))=0\,.
\end{equation}
Hence  this ideal sheaf is a cohomologically trivial line bundle.

Now we shall use \cite{BPV}, VIII, Lemmas 16.4, 17.1, 17.2, and 17.3.
These results imply that there always exist ``half pencils'' on an 
Enriques surface $X$, that is, effective divisors $D$ such that the normal 
bundle $\O_D(D)$ is a line bundle of order $2$, where $D$ is a
non-multiple divisor which
is an ``elliptic configuration'', i.e. a member of the Kodaira list in
Table 3 on p. 150 in \cite{BPV}, V, Sect.~7.  The discussion of case (c)
in the same section (cf. top of p. 151) says, in fact, that $D$ is reduced,
and is either a smooth elliptic curve, an irreducible curve with one ordinary
double point, or a polygon of $b$ curves, each of which is a smooth rational
$-2$ curve (Kodaira's types $I_0$, $I_1$ and $I_b$, respectively).

If we have a half pencil of type $I_b$ with $b\ge 2$, we have a smooth
$-2$ curve $E$ on the surface. Then as seen above, $\O_X(-E)$ is 
cohomologically trivial.

The same argument also takes care of the case of ``special Enriques 
surfaces'', on any of which there is a smooth $-2$ curve.

So we may assume our Enriques surface is ``non-special'', with two
distinct half pencils $D_1$, $D_2$, each of which is irreducible, and is
either an elliptic curve, or a singular rational curve with one ordinary
double point; further, we may assume their intersection number
$(D_1\cdot D_2)$ is $1$ (cf. \cite{BPV}, VIII, Theorem 17.7).
Thus $D_1$ and $D_2$ must intersect transversally at a point, say $x\in X$,
which is a smooth point of each curve. 

In this situation, we claim that  $\O_X(D_1-D_2)$ is a cohomologically 
trivial line bundle.

First, note that $E=D_1-D_2$ has self-intersection $-2$, since
$D_1^2=D_2^2=0$, while $(D_1\cdot D_2)=1$. Since the canonical bundle on
an Enriques surface is numerically trivial (it is $2$-torsion), and
$\chi(X, \O_X)=1$, the Riemann-Roch theorem gives that $\chi(X, 
\O_X(E))=0$.

Next, consider the exact sequence of sheaves
\begin{equation}
0\to \O_X(D_1-D_2) \to \O_X(D_1) \to \O_{D_2}(x)\to 0\,,
\end{equation}
where we note that the restriction of $\O_X(D_1)$ to the curve $D_2$ is
the line bundle corresponding to the intersection point $x$. Since $D_2$
is an irreducible curve of arithmetic genus $1$, $\O_{D_2}(x)$ has a
1-dimensional space of sections, and vanishing higher cohomology. Hence
$$
H^0(X,\O_X(D_1)) \to H^0(X,\O_{D_2}(x))\,,
$$
being a nonzero map between 1-dimensional vector spaces, is an isomorphism.
This means that $\O_X(D_1-D_2)$ has no nontrivial global sections, and
$$
\O_X(D_1-D_2)\to \O_X(D_1)
$$
is an isomorphism on higher cohomology.

Now consider the sequence
\begin{equation}
0\to \O_X \to \O_X(D_1)\to \O_{D_1}(D_1) \to 0\,,
\end{equation}
where $\O_{D_1}(D_1)$, the normal bundle to $D_1$, is a line bundle of
order 2. In particular, it a nontrivial line bundle of degree zero on a
curve of arithmetic genus $1$, and thus has no cohomology (it has no
sections, and by the Riemann-Roch theorem, its Euler characteristic
is equal to zero).
Hence $\O_X\to \O_X(D_1)$ induces isomorphisms on all cohomology groups.
In particular, $\O_X(D_1)$ has vanishing cohomology $H^1$ and $H^2$.
Hence so does $\O_X(D_1-D_2)$. 
\end{proof}

\subsection{Elliptic surfaces}
We now pass to elliptic surfaces. Here we have a result for 
elliptic surfaces with a section (i.e., for Jacobian fibrations). 

\begin{proposition}
If $f:X\to C$ is an elliptic fibration over a smooth projective curve $C$,
such that $f$ has a section, then $X$ supports a cohomologically trivial 
line bundle, and so (D) holds.
\end{proposition}
\begin{proof}
Let $D$ be a curve in $X$ mapping isomorphically to $C$ under
$f$. 

We may replace $f$ by a suitable relatively minimal model (the
pull-back of a cohomologically trivial line bundle from the minimal model
is of course cohomologicaly trivial). Since $f$ is an elliptic fibration,
all fibres of $f$ are connected curves of arithmetic genus $1$. Hence
for {\em any} fiber $F$, $H^1(F,\O_F(D))=0$, while $H^0(X,\O_X(F))$ is
1-dimensional for general fibers $F$ (note that the general fiber
is smooth). 

It follows that $R^1f_*\O_X(D)=0$ (see \cite{Ha} III Ex. 11.8, for
example). Hence for any line bundle $\L$ on $C$, the projection formula
identifies
$$
R^if_*{f^*\L\otimes_{\O_X} \O_X(D)}
$$
with
$$
\L\otimes_{\O_C} R^if_*\O_X(D)\,.
$$
Hence we obtain an isomorphism 
\begin{equation}
H^i(C,\L\otimes f_*\O_X(D)) \cong
H^i(X,f^*\L\otimes \O_X(D))
\end{equation}
for any $i\geq 0$. The sheaf $f_*\O_X(D)$ on $C$ is torsion-free, hence it 
is a vector bundle. By looking at the generic fiber, we see that it is, in 
fact, a line bundle.

We now claim that, for any line bundle $\M$ on $C$, there exists a line
bundle $\L$ on $C$ so that $\L\otimes \M$ on $C$ is cohomologically
trivial.  Indeed, by Lemma~\ref{ctl-curves}, $C$  supports a 
cohomologically  trivial line bundle $\L_0$. Choose  
$\L=\L_0\tensor\M^{-1}$. 
\end{proof}

\begin{remark}
Proposition~\ref{bielliptic} 
for $X=(E\times C)/G$, with $C$ of genus $\geq 2$, also gives examples of 
complex elliptic surfaces for which (D) holds; some of these surfaces have 
elliptic fibrations without a section. 
\end{remark}

\subsection{Some surfaces for which (D) fails}
We shall give now more surfaces for which (D) fails, because the surface 
does not support a cohomologically trivial line bundle. 

The result below applies to any surface $X$, which is a sufficiently 
general complete intersection in projective space (or in a homogeneous 
space $G/P$, where $G$ is a semisimple group, and $P$ a maximal parabolic 
sugroup), of large enough multi-degree that the canonical bundle of the
surface has non-zero sections (this amounts to saying that either $X$ is 
a K3 surface, or $X$ is of general type). This is because, since the 
surface $X$ is a very general member in the corresponding family of 
surfaces, 
we have $\Pic(X)\cong\Z$, by the Noether-Lefschetz theorem.

\begin{proposition}\label{Tz} \
Suppose $X\not\cong \bP^2$ is a smooth projective surface,
with $\Pic(X)={\Z}$, such that the ample generator of $\Pic(X)$ has a 
non-zero section. Then $X$ does not have (D). 
\end{proposition}
\begin{proof}
Note first that (\ref{prod}) holds by Lemma \ref{Pic}. So we must show 
that $X$ does not support any  cohomologically trivial line bundle, that 
is, every power of the ample generator has some non-zero cohomology.

If $\L$ is a nonnegative power of the ample generator, then clearly 
$H^0(X,\L)\ne 0$, by our assumptions. Suppose now that $\L$ is a negative 
power of the ample generator. Observe that $\omega_X$ is a nonnegative 
power of the ample generator because (by classification) $X$ is not Fano, 
since the only Fano surface with Picard group $\Z$ is $\bP^2$. We conclude
that $\L^{-1}\otimes \omega_X$ is a strictly positive power of the ample 
generator, and so, by Serre duality,
$$
H^2(X,\L)=H^0(X,\L^{-1}\otimes \omega_X)\ne 0\,.
$$
The proposition has been proved.
\end{proof}

\noindent
\begin{remark}
We comment on the hypothesis in Proposition~\ref{Tz} that the ample
generator of $\Pic(X)$ has a non-zero section. Any surface 
$X\not\cong\bP^2$ with $\Pic(X)=\Z$ must have trivial canonical bundle, or
be of general type, with ample canonical bundle $\omega_X\cong \O_X(n)$ 
for some $n>0$, where $\O_X(1)$ denotes the ample generator of $\Pic(X)$.

We claim that the canonical bundle of such a complex surface always has a 
non-zero section, i.e. $p_g(X)>0$.

Indeed, suppose that the canonical bundle is ample, with no sections.
It follows from Hodge theory that the second Betti number of $X$ equals
its Picard number. Hence the topological Euler characteristic of $X$ is
$3$, while its holomorphic Euler characteristic is $1$. This forces 
$K_X^2=9$ by the Noether formula, that is, $X$ is a {\em Fake $\bP^2$} in 
the sense of  Mumford \cite{M1}, who first constructed such a surface of 
general type, using $p$-adic uniformization.

All complex Fake $\bP^2$'s have been classified in recent work of Gopal
Prasad and Sai-Kee Yeung \cite{PY}. From their work (see Theorem~10.1),
every Fake $\bP^2$ has a non-zero torsion subgroup of $\Pic(X)$, and
hence {\em does not} have Picard group $\Z$.

It is interesting to ask if a complex Fake $\bP^2$ has (D). Poincar\'e 
duality, combined with the Riemann-Roch theorem, gives us that any ample
generator of the N\'eron-Severi group, modulo torsion, has
self-intersection 1, and thus vanishing Euler characteristic. Hence (D)
holds for $X$ precisely when a line bundle can be found on $X$
which gives an ample generator for $NS(X)/({\rm torsion})$, and has 
vanishing $H^1$.
\end{remark}

\begin{remark}
We do not know an example of a surface $X$ with $\Pic(X)=\Z$, and an ample
canonical bundle with a non-zero section, but with $H^0(X,\O_X(1))=0$ for
the ample generator $\O_X(1)$ of $\Pic(X)$.
\end{remark}

\section{Higher dimensional varieties}\label{ptpr}

\subsection{Varieties with Picard group $\Z$}

We first consider varieties of dimension $d\geq 3$ with Picard group $\Z$.
From the Grothendieck-Lefschetz theorem (cf., e.g., \cite{Har1}), we
have $\Pic(X)=\Z$ for any smooth complete intersection $X$, or for a
smooth complete intersection of divisors in a homogeneous space $G/P$
with $G$ semisimple, $P$ maximal parabolic (so that $\Pic(G/P)=\Z$).

\begin{proposition}\label{fano1}
Let $X$ be a smooth projective variety of dimension $d\geq 3$ over a
field. Suppose that
\begin{enumerate}
\item[(i)] (D) holds for $X$;
\item[(ii)] $\Pic(X)=\Z$, and the ample generator of $\Pic(X)$ has a
nonzero section.
\end{enumerate}
Then $X$ is a Fano variety, with canonical line bundle $\omega_X\cong
\O_X(-n)$ for some $n\geq 2$.
\end{proposition}
\begin{proof}
Since $X$ has Picard group $\Z$, and (D) holds, Theorem~\ref{Tm} and
Lemma~\ref{Pic} imply that $X$ supports a cohomologically trivial line
bundle, which must be $\O_X(m)$ for some integer $m$. If
$\omega_X=\O_X(r)$  then by Serre duality, $\O_X(m)$ and $\O_X(r-m)$ both
have $H^0=0$. Since the ample generator of $\Pic(X)$ has a nonzero
section, we must have  $m<0$, $r-m<0$. Hence $r=(r-m)+m\leq -2$.
\end{proof}
\begin{corollary} Let $X\subset \bP^n$ be a smooth complete intersection
of multidegree $(d_1,\ldots,d_r)$ with $r\leq n-3$, and $\sum_id_i\geq
n$. Then $X$ does not have (D).
\end{corollary}
\begin{proof} Such a complete intersection has dimension $n-r\geq 3$, so
has Picard group $\Z$, from the Grothendieck-Lefschetz theorem. The
canonical bundle of $X$ is
$$
\O_X\Bigl(\sum_id_i-n-1\Bigr)\,,
$$
where $\sum_id_i-n-1\geq -1$. Hence (D) does not hold for $X$.
\end{proof}
Thus, for smooth hypersurfaces in $\bP^4$, the only cases we need to
consider are quadrics and cubics; the case of quadrics is settled below
((D) does not hold, cf. Proposition~\ref{quadric1}).

For Fano varieties, we similarly have:
\begin{corollary} Let $X$ be a smooth projective complex Fano variety with
second Betti number $b_2(X)=1$ and $\omega_X=\O_X(-1)$ (i.e. $X$ is of
index 1). Then $X$ does not have (D).
\end{corollary}
\begin{proof}
Recall that $H^i(X,\O_X)=0$ for $i>0$, by Serre duality and the
Kodaira vanishing theorem; in particular, $X$ is algebraically simply
connected, since $\chi(X,\O_X)=1$ . Since $b_2(X)=1$, we have
$\Pic(X)={\Z}$, and (\ref{prod}) holds (cf. Lemma
\ref{Pic}). Now Proposition~\ref{fano1} implies that (D) does not hold.
\end{proof}
Fano varieties of the above type have been essentially classified by
Iskovskih (see \cite{Is}, \cite{Cu}). This is extended to positive
characteristics in \cite{SB}.

\subsection{(D) and the Point Property}
It is useful to consider a property related to (D), which is sometimes a
consequence of it.

Let $X$ be a scheme. For a line bundle $\L$ on $X$, we say that the
``$\L$--point property holds'' if the following is true:

\smallskip

`` If for each $x\in X$, there exists a vector bundle $\F$ on $X$ of rank
$d=\dim(X)$ with $\det(\F)=\L$, and a section of $\F$ vanishing exactly
at $x$ with multiplicity $1$.''

\medskip

\noindent
This implies
$$
c_1(\F)=c_1(\L)\in CH^1(X),\;\;c_d(\F)=[x]\in CH^d(X) \,.
$$

\begin{theorem}
Let $X$ be smooth and proper over and algebraically closed field $k$.
Denote by $p_1, p_2: X \times X \to X$ the two projections. Suppose that 
the diagonal property (D) holds, and $\Pic(X)$ is finitely generated. Then 
there exists a cohomologically trivial line bundle $\L$ on $X$ such that

\smallskip

\noindent
(i) the $\L^{-1}$--point property holds, and also

\smallskip
\noindent
(ii) the $\L\otimes \omega_X^{-1}$--point property holds.
\end{theorem}
\begin{proof}
Let $\dim(X)=d$, and let $\E$ be a rank $d$ bundle on $X \times X$
given by (D).

By lemma~\ref{Pic}, the finite generation of $\Pic(X)$ is equivalent to
(\ref{prod}).  By (\ref{prod}), we have
\begin{equation}\label{produ}
\det(\E)= p_1^*\L_1 \otimes p_2^*\L_2\,,
\end{equation}
for some line bundles $\L_1$ and $\L_2$ on $X$. Then the restrictions
of $\E$ to $X\times\{x\}$ and $\{x\}\times X$ have the determinants
$\L_1$, $\L_2$ respectively, whatever the choice of $x$. Hence the
$\L_1$--point property as well as the $\L_2$--point property hold for $X$.

When the determinant of $\E$ does have the above special form
(\ref{produ}), then the restriction to the diagonal is
$\L_1\otimes \L_2$, so that, by Theorem~\ref{Tm}, we see that
$\L_1^{-1}$ and $\L_2^{-1}$ are mutually Serre dual line
bundles, which are cohomologically trivial line bundles.
\end{proof}
\smallskip

\begin{corollary}\label{Lpt} Assume that a smooth proper $k$-scheme
$X$ satisfies (\ref{prod}). If for any cohomologically trivial line
bundle $\L$ on $X$, either the $\L^{-1}$--point property fails, or the
$\L\otimes \omega_X^{-1}$--point property fails, then $X$ does not admit (D).
\end{corollary}

Observe, that there exists a ``trivial'' instance of the corollary
when there are {\it no} cohomologically trivial line bundles on $X$ at all.

\begin{proposition}\label{quadric1} A smooth quadric $Q_3\subset {\Bbb
P}^4$
over an algebraically closed field fails to have (D).
\end{proposition}
\begin{proof}
There are two cohomologically trivial line bundles on $Q_3$: $\L_1=\O_{Q_3}(-1)$
and $\L_2=\O_{Q_3}(-2)$. Since $\omega_{Q_3}=\O_{Q_3}(-3)$, we have
$$
\L_1^{-1}=\L_2 \otimes \omega_{Q_3}^{-1}\,,
$$
so by Corollary~\ref{Lpt}, it suffices to show that
the $\L_1^{-1}$--point property fails.

We use a standard presentation
of the Chow ring $CH^*(Q_3)$ as ${\Z}^{\oplus 4}$ with a suitable ring
structure. Let $[Q_2]$, $[L]$, and $[P]$ (quadric surface, line, and point)
be the generators of $CH^1(Q_3)$, $CH^2(Q_3)$, and $CH^3(Q_3)$ respectively.
There are the following relationships:
\begin{equation}
[Q_2]^2=2[L] \ \ \ \hbox{and} \ \ \ [Q_2]\cdot [L]=[P]\,.
\end{equation}

If $\E$ is a vector bundle on $Q_3$, the total Chern class of $\E$
is of the form
$$
1+d_1(\E)[Q_2]+d_2(\E)[L]+d_3(\E)[P]\,,
$$
where $d_i(\E)\in {\Z}$.

The argument now boils down to showing that there is no rank $3$ vector
bundle $\E$ on $Q_3$ with $d_3(\E)=1$ and $d_1(\E)=1$.

Suppose - ad absurdum - that such a bundle exists. We use the formula for
the Euler characteristic of $\E$ given by the
Grothendieck-Hirzebruch-Riemann-Roch theorem (see \cite{F2},
Example~15.2.5 for a general formula for 3-folds).
In fact, we use the following explicit version of the formula for smooth
quadric 3-folds, given in \cite{ES}: for a vector bundle $\E$ on $Q_3$,
\begin{equation}\label{hrr}
\chi(Q_3,\E)=\frac{1}{6}(2d_1^3-3d_1d_2+3d_3)+\frac{3}{2}(d_1^2-d_2)
+\frac{13}{6}d_1 +\rank(\E)\,,
\end{equation}
where $d_i=d_i(\E)$ are the above numbers. Substituting the present
values, we get
$$
\chi(Q_3,\E)=\frac{15}{2}-2d_2 \,.
$$
This contradicts the fact that $\chi(Q_3,\E)$ is integer, and
the proposition has been proved.
\end{proof}

Note that since the quadric $Q_3$ is birationally isomorphic to the
projective $3$-space which has (D), this last property is not a birational
invariant in dimension $\ge 3$.

\smallskip

The following result is obtained in a standard way, following the
argument that the tangent bundle of a group variety is trivial.

\begin{proposition}\label{Pwpp} Let $X$ be a group variety over an algebraically
closed field. Then $X$ has (D) if and only if $X$ has the following
``weak point property'': for some point $x\in X$,
there exists a vector bundle $\E$ of rank $d=\dim(X)$, such that there is a
section of $\E$ with zero scheme $x$.
\end{proposition}
\begin{proof} If $X$ has (D), then for any $x\in X$, the weak point
property holds, by restriction of the data giving (D) to $X\times \{x\}$.

Conversely, suppose the weak point property holds with respect to the
point $x$. Let $\E$ be the corresponding vector bundle, and $s$ a section
with zero scheme $x$. Let $\mu:X\times X\to X$ be the multiplication, and
$i:X\to X$ the inverse, defining the algebraic group structure on $X$.
Define a new map $f:X\times X\to X$ by
\[f(u,v)=\mu(\mu(u,i(v)),x).\]
Then $f$ is a morphism, whose scheme theoretic fiber
$f^{-1}(x)=\Delta_X$, the diagonal subscheme. Hence the vector
bundle $f^*\E$ has the section $f^*(s)$ whose zero scheme is the diagonal.
\end{proof}
\begin{corollary} An abelian variety has (D) precisely if it has the weak
point property for some point.
\end{corollary}

Though we have not resolved whether a smooth cubic 3-fold has (D), we have
the following observation, which is a necessary condition for (D), since
the only cohomologically trivial line bundle on a cubic 3-fold $X$ is
$\O_X(-1)$.
\begin{lemma} A smooth complex cubic 3-fold $X$ has the $\O_X(1)$-point
property.
\end{lemma}
\begin{proof} It is known that a smooth cubic 3-fold $X$ contains lines
from the ambient projective space, and is covered by such lines. Further,
any line $L\subset X$ has a normal bundle with trivial determinant, from
the adjunction formula, since $\omega_X=\O_X(-2)$, $\omega_L=\O_L(-2)$,
and $\O_X(1)\tensor \O_L=\O_L(1)$ (the last formula holds because $L$ is a
line).

If $\I_L$ is the ideal sheaf of a line in $X$, then since $H^i(X,\O_X)=0$,
$i=1,2$, we have an isomorphism
\[{\Ext}^1(\I_L,\O_X)\cong
H^0(\underline{\Ext}^1(\I_L,\O_X))\cong\Hom(\det\I_L/\I_L^2,\O_L)=\C.\]
Thus, from the Serre argument,  there is a rank 2 vector bundle $\F$ on
$X$ with trivial determinant, together with a section, whose zero scheme
is $L$.

Now suppose $x\in X$ is any point. Choose a line $L\subset X$ passing
through $x$, and choose a hyperplane $H$ in $\bP^4$, not containing $L$,
but passing through $x$. Then $L$ and $H$ intersect transversally at $x$.
If $t\in H^0(\bP^4,\O(1))$ is a section with zero scheme $H$, then the 
rank 
3 vector bundle $\E=\F\oplus\O_X(1)$ has a section $(s,t\mid_X)$ with zero 
scheme $L\cap(H\cap X)=\{x\}$, and the intersection is transverse in $X$.
The determinant of $\E^*$ is clearly $\O_X(-1)$.
\end{proof}

\section{The affine case}
It makes sense to ask which smooth affine varieties over an algebraically
closed field have (D). If $X=\Spec A$ is a smooth affine $k$-variety of
dimension $\leq 2$, then $X$ has (D); in the 2-dimensional case, this
follows from Serre's classic argument, since the diagonal is a smooth,
codimension 2 subscheme of an affine scheme.

We do have one positive result in the affine case, which may be of
interest. This is basically a corollary of work of M. P. Murthy.
\begin{proposition}\label{wpp} \ An affine algebraic group over an algebraically
closed field has (D).
\end{proposition}
\begin{proof}
If  $X$ is an affine algebraic group over
an algebraically closed field, then a result of M. P. Murthy (cf.
\cite{Mur}, Theorem 3.3) shows that any smooth point is the zero scheme of
a section of a vector bundle of rank equal to $\dim(X)$.
\end{proof}
\begin{remark}
In fact for an affine algebraic group $X$, we can even say that any point
$x\in X$ is a complete intersection (i.e. the vector bundle may be chosen
to be a {\it trivial bundle}). Hence the diagonal is also a (global)
complete intersection. Indeed, Murthy's results \cite{Mur} imply that a
point on a smooth affine variety is a complete intersection if and only if
its class in the Chow group of points vanishes. This holds for $X$
because for unirational smooth affine varieties the Chow group of points
is, in fact, zero. One can now naturally ask if the diagonal of any affine
algebraic group is a complete intersection in any group-theoretically
significant way!
\end{remark}

We conjecture that there exist smooth affine complex varieties of
any dimension $\geq 3$ for which (D) fails. However, we have been unable
to construct such examples. This leads to the following question:\\[2mm]
{\bf Question}:\quad Let $A$ be a smooth affine algebra over an
algebraically closed field $k$. Let $K$ be an extension field of $k$ (not
necessarily algebraically closed), $A_K=A\tensor_kK$, and and $M\subset
A_K$ a maximal ideal with residue field $K$. Does there exist a projective
$A_K$-module $P$ of rank $n=\dim(A)$ such that there is a surjection $P\to
M$?

\bigskip

If $K$ is an algebraically closed extension of $k$, this is always true,
from Murthy's results. If $X=\Spec A$ satisfies (D), then the question
has a positive answer for any field extension $K$. So a negative answer to
the question would give a way  of constructing counterexamples to (D).

\section{The Topological Diagonal Property}
Let $M$ be a smooth compact connected oriented manifold of real dimension $n$\,,
and let
$$
\Delta\subset M\times M
$$
be the diagonal submanifold. In this section, we introduce
topological versions of the diagonal property $(D)$. We also use
notation standard in topology (instead of that from algebraic geometry
used earlier).

We say that $M$ {\it has property} $(D_{r})$ if there exists
a smooth real vector bundle $\E$ of rank $n$ on $M\times M$ and
a smooth section
$s$ of $\E$ such that (i) $s$ is transverse to the $0$-section $0_{\E}$ of
$\E$ and (ii) $\Delta = s^{-1}(0_{\E})$. If further the bundle $\E$
is {\em orientable}, we say that $M$ {\it has property} $(D_{o})$. Finally, if
$\dim_{\R}M=2m$, and $\E$ can be chosen to be a smooth complex vector bundle
of ${\rank}_{\C}\E=m$, we say that $M$ {\it has property} $(D_{c})$.

\begin{remark}{\rm Clearly $(D_{o})\Rightarrow (D_{r})$. Further, if $M$ satisfies
$(D_{c})$, then $\E_{|\Delta}$ is isomorphic to the
normal bundle of $\Delta$ by the transversality condition,
and this normal bundle is known to be isomorphic to the tangent bundle
$\tau_{M}$ (cf. (\ref{con}) in \S2). It follows that $\E$ a complex
vector bundle forces $\tau_{M}$ to be a complex complex bundle, so $M$
will have to be an {\em almost complex manifold}. $(D_{o})$ and $(D)$ will
follow
from $(D_{c})$ whenever it holds. Also if $M$ is a complex manifold,
then $(D)$ of the earlier sections will imply $(D_{c})$, hence $(D_{o})$
and $(D_{r})$.}
\label{almcomp}
\end{remark}

\begin{remark}(\rm The topological point property)
{\rm In analogy with ``weak point property'' from Section 4, if $M$ has $(D_{r})$,
then letting $\E$ be the bundle realizing $(D_{r})$, and setting
$\E_{1}:=\E_{M\times\{p\}}$, and the section $\sigma =
s_{|M\times\{p\}}$ one finds
$\sigma\transv 0_{\E_{1}}$ and $\sigma^{-1}(0_{\E_{1}})=\{p\}$. Likewise for
$\E_{2}:=\E_{|\{q\}\times M}$. We will say that a smooth manifold {\em obeys $(P_{r})$ (resp.
$(P_{o})$ resp. $(P_{c})$) if there exists a smooth real (resp. real orientable, resp. complex)
bundle $\H$ of ${\rank}_{\R}=\dim_{\R}M=n$ (resp.
${\rank}_{\R}=\dim_{\R}M=n$,
resp. ${\rank}_{\C}\H=n$) with a smooth section $\sigma$
meeting $0_{\H}$ transversely at one point}. Thus $(D_{r})$ (resp. $(D_{o})$, resp.
$(D_{c})$) implies $(P_{r})$ (resp. $(P_{o})$, resp. $(P_{c})$). If $M$ is compact,
by well known transversality results, $\H$ realizes $(P_{r})$ (resp. $(P_{o})$, resp.
$(P_{c})$) iff $w_{n}(H)=1$ (resp. $e(\H)=\pm 1$, resp. $c_{n}(\H)=\pm 1$).}
\label{euler1}
\end{remark}

\begin{example}{\rm (Lie Groups) Let $G$ be a Lie group, then $G$ satisfies
$(D_{r})$ (resp. $(D_{o})$, resp. $(D_{c})$) iff it satisfies $(P_{r})$ (resp.
$(P_{o})$, resp. $(P_{c})$.) The proof is analogous to Proposition \ref{wpp}.
Let $\H$ be a bundle realizing $(P_{r})$ and a section with a transverse zero
at $1\in G$ (by translating if necessary), and then the bundle
$\E:=\mu^{\ast}\H$\,, where $\mu:G\times
G\rarr G$ is $(x,y)\mapsto xy^{-1}$ with pulled back section realizes $(D_{r})$. Likewise for
$(D_{o})$ and (if $D\mu$ is $\C$-linear) also $(D_{c})$. It follows that $S^{1}$ satisfies $(D_{r})$.
Similarly $SO(3)=\R\bP^{3}$ satisfies $(D_{r})$. We will
see later that no odd dimensional manifold satisfies $(D_{o})$, so
neither of the above satisfies $(D_{o})$.}
\label{liegroup}
\end{example}

\begin{example}{\rm (Products) If two manifolds $M$ and $N$ have
$(D_{r})$ (resp. $(D_{o})$, resp. $(D_{c})$, then so does $M\times N$,
as remarked in the Introduction, by taking the product of the corresponding
bundles. In particular, by the example of $S^{1}$ in Example \ref{liegroup} above, and
and of $\C\bP^{1}$ from the next example, manifolds satisfying
$(D_{r})$ (resp. $(D_{c})\Rightarrow (D_{o})$) exist in every dimension
(resp. every even dimension).}
\label{prodexample}
\end{example}

\begin{example}{\rm (Projective spaces and Grassmannians) For all $k$
the real Grassmannian $G_{k}(\R^{n})$ satisfies $(D_{r})$, by the
construction using the tautological bundles given in the Introduction. Likewise $G_{k}(\C^{n})$
satisfies $(D_{c})$ by the analogous complex construction. It is not
true that an orientable Grassmannian satisfies $(D_{o})$, e.g. we
shall see below that no odd-dimensional real projective space
satisfies $(D_{o})$.}
\label{projgras}
\end{example}

\begin{example}{\rm (Compact Riemann Surfaces) As remarked in the Introduction, if
$M$ is a compact connected
Riemann surface, then it satisfies $(D)$, and hence $(D_{c})$.}
\label{riemsurf}
\end{example}

\begin{example}{\rm (The spheres $S^{1}, S^{2}, S^{4}, S^{8}$) The
circle $S^{1}$, in view of Example \ref{liegroup} (or \ref{projgras})
satisfies $(D_{r})$, and $S^{2}=\C\bP^{1}$ satsifies $(D_{c})\Rightarrow
(D_{o})$ in view of Example \ref{projgras}. Regarding $S^{4}$ (resp.
$S^{8}$) as the quaternionic line (resp. octonionic line) and
repeating the construction of Example \ref{projgras} above with respect to the
respective tautological bundles defined on these, one checks that $S^{4}$ and
$S^{8}$ both satisfy $(D_{o})$. They obviously can't satisfy $(D_{c})$
because neither is an almost complex manifold (cf. Remark \ref{almcomp}).}
\label{projline2}
\end{example}

A Riemannian metric on $M$ induces  one on $M\times M$, so metrics result
on $\tau_{M\times M}$, and all its subbundles. Let $U$ be a closed $\eps$-tubular
neighborhood of $\Delta$ in $M\times M$. The tubular neighborhood theorem
gives a smooth diffeomorphism $\phi: (U,\boun U)\rarr (D(\nu),S(\nu))$, where
$D(\nu)$ is the $\eps$-disc bundle of the normal bundle $\rho:\nu\rarr \Delta$ of $\Delta$ in
$M\times M$, and $S(\nu)$ the $\eps$-sphere bundle. The map $r:=\rho\circ \phi :U\rarr \Delta$ 
is then a strong deformation retraction of $U$ to its core $\Delta$. Thus there is a bundle diagram
\begin{eqnarray}
\begin{array}{ccc}
r^{\ast}(\nu)&\stackrel{D\phi}{\longrightarrow}& \rho^{\ast}(\nu)\\
\rho\downarrow & &\downarrow\rho\\
U&\stackrel{\phi}{\longrightarrow}&D(\nu)\,.
\end{array}
\label{diag1}
\end{eqnarray}
The restricted bundle $\rho^{\ast}(\nu)_{|S(\nu)}\rarr S(\nu)$ has a tautological
section $s$ defined by
$v\mapsto v$, which satisfies $\norml s(v)\normr=\eps$ for all $v\in
S(\nu)$. Thus we have an orthogonal direct sum decomposition of bundles 
on $S(\nu)$:
$$
\rho^{\ast}(\nu)_{|S(\nu)}= \xi\oplus \L\,,
$$
where $\L$ is the trivial line subbundle spanned by $s$ (and $\xi$ denotes
its orthogonal complement). It is also well known that $\rho:\nu\rarr \Delta$ is
isomorphic to the tangent bundle $\rho:\tau_{M}\rarr M$ under the identification
$\Delta\cong M$, and since $M$ is orientable, so is $\xi$. Hence, under
this identification,
$\xi$ is isomorphic to the quotient bundle
$$
\rho^{\ast}(\tau_{M})/\L\rarr S(\tau_{M})\,,
$$
where $\L\rarr S(\tau_{M})$ is the trivial tautological bundle spanned by the tautological
section of $\rho^{\ast}(\tau_{M})$ over $S(\tau_{M})$.

Let $\F:=\phi^{\ast}(\xi)$, a rank $(n-1)$ subbundle of $r^{\ast}(\nu)_{|\boun U}$.
It is isomorphic to the rank $(n-1)$ bundle $\rho^{\ast}(\tau_{M})/\L\rarr
S(\tau_{M})$ under the above identifications. Note that $\F$ is an
orientable bundle on $\boun U$.

\begin{remark}{\rm The restriction of the bundle $\xi$ above to each fiber
$S(\nu_{x})$ of the sphere bundle $\rho:S(\nu)\rarr M$
is the tangent bundle $\tau_{n-1}$ of the sphere $S(\nu_{x})$. Consequently, the bundle
$\F$, when restricted to a fiber $r^{-1}(x)$ of the fiber bundle $r:\boun
U\rarr \Delta$, is isomorphic to $\tau_{n-1}$. This is clear, since the
fiber of $\L$ at a point $v\in S(\nu_{x})$ is precisely $\L_{v}=\R v$}.
\label{remark1}
\end{remark}

\smallskip

Since the $n=1$ (i.e. $S^{1}$) case is completely settled by Example
\ref{liegroup}, we will assume henceforth that $n=\dim_{\R}M\geq 2$. The following
lemma is the {\it key technical} result of this section. We shall write
$U^{\circ}$ for the interior of $U$.

\begin{lemma}\label{lemma4} Let $M$, $U$, $\Delta$, be as above. Set
$X:=(M\times M)\setminus U^{\circ}$.
Then $M$ has $(D_{r})$ iff
the rank $(n-1)$ bundle $\F\rarr \boun U$ defined above is isomorphic
to the restriction to $\boun U=\boun X$ of a smooth rank $(n-1)$ bundle
$\G$ on $X$. Further, $M$ has $(D_{o})$ iff the bundle $\G$ can be chosen to
be orientable.

These are both problems in homotopy theory, of extending
classifying maps $\boun X\rarr BG$ for the bundle $\F$ (where $G=O(n-1)$ for $(D_{r})$,
resp. $SO(n-1)$ for $(D_{o})$) to $X$.
\label{lemma1}
\end{lemma}
\begin{proof}
We first prove the ``only if'' part. Suppose there exists a
rank $n$ smooth real vector bundle $\pi:\E\rarr M\times M$, and $s$ a
smooth section transverse to the zero-section $0_{\E}$ such that the
diagonal $\Delta = s^{-1}(0_{\E})$. The strong deformation retraction
$r:U\rarr \Delta$ makes $\E_{|U}$
isomorphic to $r^{\ast}(\E_{|\Delta})$, so $\E_{r(y)}$ is identified with $\E_{y}$ for all
$y\in U$. The section $s$ is
nowhere vanishing on
$$
X=(M\times M)\setminus U^{\circ}\,,
$$
and hence defines a trivial line subbundle $\Lambda$ of $\E_{|X}$, and a splitting of bundles on $X$:
$$
\E_{|X}= \G\oplus \Lambda\,,
$$
where $\G$ is a rank $(n-1)$ bundle on $X$.
Using the transversality condition on $s$, the tubular neighborhood theorem, the
inverse function and compactness
of $M$, and choosing $\eps$ (the radius of $U$) small enough, it is easy to see that there is a
smooth bundle
equivalence of $\pi:\E_{|U}\rarr U$ and $\rho:\tau_{\Delta|U}\rarr U$ which carries the
section $s$ of the former to the tautological section of the latter ($U$ being identified with
the $\eps$-disc bundle $D(\tau_{\Delta})=D(\tau_{M})=D(\nu)$). Thus on $\boun U=\boun X$,
quotienting by the respective trivial line bundles defined by these sections, it
follows that $\G_{|\boun
U}=\G_{|\boun X}$ is isomorphic to  $\F$. If $M$ had $(D_{o})$, and $\E$ was
an orientable bundle realizing $(D_{o})$
to begin with, then $\G=\E_{|X}/\eps^{1}$ would also be orientable.
This proves the ``only if'' parts of both the first and second statements.

For the ``if'' part, let $\G$ be given on $X$ as in the statement. Construct the bundle $\E$,
by taking $r^{\ast}(\nu)$ on $U$, and gluing it
to the bundle $\G\oplus \eps^{1}_{X}$ on $X$, after ensuring that the
decomposition $r^{\ast}(\nu)=\F\oplus \phi^{\ast}(\L)$ on $\boun U$ is preserved, viz. the first summand
$\G_{|\boun U}$ is glued
to the first summand $\F_{|\boun U}$ via the given isomorphism, and the second trivial
summand $\eps^{1}_{X}$ is
glued to the second trivial summand
$\phi^{\ast}\L_{|\boun U}$ on $\boun U$ by matching
the section $\sigma$ of the trivial bundle $\phi^{\ast}\L\rarr \boun U$ defined above
with any everywhere $\neq 0$
section of $\eps^{1}_{X}$ which extends $\sigma$ to $X$ (this is possible by
the theorem of Tietze, since $n\geq 2$ implies $\boun X$ is connected).
The section $\sigma$ of $r^{\ast}(\nu)\rarr \boun U$ above is the
restriction of the tautological section, also denoted  $\sigma$,
of $r^{\ast}(\nu)\rarr U$ which is transverse to the zero section.
Thus the matched section $s$ of the whole bundle $\E$ vanishes exactly on
$\Delta$, with $s\transv 0_{\E}$.

The lemma has been proved.
\end{proof}

\smallskip

We have an analogue of this lemma for
$(D_{c})$. Assume (in view of Remark \ref{almcomp}) that $M$ is an almost complex manifold of
real dimension $n=2m$.
Then in the notation of the last section, the normal bundle $\nu\cong \tau_{M}$ of $\Delta$
is a complex vector bundle, and $r^{\ast}(\nu)\rarr \boun U$ splits off a {\em complex} line subbundle
$\eps^{1}_{c}$ defined by the complex span of the tautological section
of $r^{\ast}(\nu)_{|\boun U}$.
Thus we may write
$$
r^{\ast}(\nu)_{|\boun U}=\F_{c}\oplus \eps^{1}_{c}\,,
$$
where $\F_{c}$ is now a complex vector bundle on $\boun U=\boun X$ with
${\rank}_{\C}\F_{c}=m-1$.

\begin{lemma} Let $M$ be an almost complex manifold of $\dim_{\C}M=m$.
Then $M$ has $(D_{c})$ iff the bundle $\F_{c}$ of complex rank
$m-1$ on $\boun X$ is isomorphic as a complex bundle to the restriction
of a complex vector bundle $\G_{c}$ on $X$. This is again a homotopy
problem as in Lemma \ref{lemma1}, with structure group $G=U(m-1)$.
\label{lemma1c}
\end{lemma}
\begin{proof}
The proof is a minor modification of that of Lemma \ref{lemma1}, and therefore is omitted.
\end{proof}

\begin{corollary}{\rm (Compact Riemann surfaces again)
Let $M$ be a compact Riemann surface. Then $M$ has $(D_{c})$}.
\label{prop2}
\end{corollary}
\begin{proof}
We saw this in Example \ref{riemsurf}. It also immediately follows
from Lemma \ref{lemma1c}, since $\F_{c}$ is a complex vector bundle of rank $0$ (!)
\end{proof}

\begin{remark}{\rm Lemmas \ref{lemma1}, \ref{lemma1c} allow us to
identify the obstructions to  $(D_{r})$, $(D_{o})$ and $(D_{c})$ as
some relative cohomology classes of
the pair $(X,\boun X)$ with local coefficients in $\pi_{i}(BG)$, where $G=O(n), SO(n), U(n)$
respectively. Since these are not computable obstructions, we skip the
details.}
\label{homotopy}
\end{remark}

\begin{theorem} The sphere $S^{n}$ has $(D_{r})$ if and only if
$n=1, 2, 4$ or $8$. (All except the first have $(D_{o})$.)
\label{prop1}
\end{theorem}
\begin{proof}
The unit sphere bundle $S(\tau_{n})$ is the Stiefel manifold
$V_{2}(\R^{n+1})$ of orthonormal $2$-frames $(x,v)$ in $\R^{n+1}$.
The bundle projection $\rho: S(\tau_{n})\rarr S^{n}$ is  projection
into the first factor, and thus we have the spherical fiber bundle
\begin{eqnarray}
S^{n-1}_{x}\stackrel{j_{x}}{\rarr}V_{2}(\R^{n+1})\stackrel{\rho}{\rarr} S^{n}
\label{stiefel}
\end{eqnarray}
with fiber $S^{n-1}_{x}=\rho^{-1}(x)$ over $x$. The tautological section of the bundle
$$
\rho^{\ast}(\tau_{n})\rarr V_{2}(\R^{n+1})
$$
is the map $(x,v)\mapsto v$. Set $\eps=\pi$ (length of a semicircle), and denote the closed
$\pi$-tubular neighborhood of $\Delta$ by $U$ as above. The complement
$$
X=(S^{n}\times S^{n})\setminus U^{\circ}
$$
is a closed $\pi$-tubular neighborhood of the antidiagonal $\Gamma$, which is the graph
of the antipodal map $A:S^{n}\rarr S^{n}$ defined by $Ax=-x$.

The involution $1\times A$ smoothly identifies $U$ with $X$.
The common boundary $\boun U=\boun X$ is diffeomorphic to the $\pi$-sphere
bundle $S_{\pi}(\tau_{n})=V_{2}(\R^{n+1})$.  Also, $X$ becomes a disc
bundle $D_{\pi}(\nu_{\Gamma})=D_{\pi}(\tau_{n})$, and contains the antidiagonal
$\Gamma$ as a strong deformation retract. Let $\theta:X\rarr \Gamma$ denote the
retraction, coming from the bundle projection 
$$
\nu_{\Gamma}\cong \tau_{n}\rarr \Gamma\,.
$$
Because $\theta$ is a deformation retraction, every bundle $\G$ on $X$ is
the $\theta$-pullback of a bundle on $\Gamma$, equivalently a $\rho$-pullback of a bundle
on $S^{n}$. Hence it follows by
Lemma \ref{lemma1} that $S^{n}$ has $(D_{r})$
iff the bundle $\F\rarr V_{2}(\R^{n+1})$ is isomorphic to the pullback under $\rho$ of some bundle on
$S^{n}$. By (\ref{stiefel}),
$$
\rho\circ j_{x}:S^{n-1}_{x}\rarr S^{n}
$$
is the constant map to $x$, it follows that the $\rho$-pullback of any bundle $\G$ on $S^{n}$
will be trivial when restricted to a fiber $S^{n-1}_{x}$. It follows
that $\F_{|S^{n-1}_{x}}$ is isomorphic to a trivial bundle on $S^{n-1}_{x}$.
By Remark \ref{remark1}, $\F_{|S^{n-1}_{x}}$ is isomorphic to the
tangent bundle $\tau_{n-1}$ of $S^{n-1}_{x}$. Hence $\F$ will be a $\rho$-pullack of a
bundle on $S^{n}$ only if $\tau_{n-1}$ is trivial. Hence $n-1=0,1,3,$ or $7$ \ by \cite{Mi}, Theorem 2.
This proves the ``only if'' part of the theorem. For the ``if'' part, cf. Example \ref{projline2}.

\smallskip

The theorem has been proved.
\end{proof}

\smallskip

All homologies and cohomologies hereafter are with $\Z$ coefficients unless otherwise stated.

\begin{theorem}
Let $M$ be a compact orientable manifold of odd dimension.
Then $M$ does not have $(D_{o})$. If further
$H_{1}(M,\Z_{2})=0$, then it does not have $(D_{r})$.
\label{prop3}
\end{theorem}
\begin{proof}
 Let $\E$ be an orientable bundle of odd real rank $2k+1$ on
$M\times M$, where $\dim(M)=2k+1$. It is well-known that the Euler
class $e(\E_{1})$ of $\E_{1}:=\E_{|M\times
\{p\}}$ must be zero, since $\E_{1}$ is of odd rank = $2k+1 =\dim(M)$ and top homology
of $M$ is $\Z$, devoid of 2-torsion. This contradicts property
$(P_{o})$ of Remark \ref{euler1}, so $M$ doesn't satisfy $(D_{o})$.

For the second assertion, note that if  $H_{1}(M,\Z_{2})=0$, then
$H^{1}(M,\Z_{2})=0$ and also $H^{1}(M\times M,\Z_{2})=0$. In particular,
every bundle on $M\times M$ is orientable. Thus if $M$ satisfies $(D_{r})$, it
automatically satisfies $(D_{o})$. But this contradicts the first statement.
\end{proof}

\begin{remark} {\rm We note from Example \ref{projgras} of
$\R\bP^{3}$ (satisfying $(D_{r})$ but not $(D_{o})$)
above that this is a
sharp result, i.e. the $H_{1}(M,\Z_{2})=0$ condition cannot be dropped. Note that
Theorem \ref{prop1} for spheres of odd dimension $n\geq 2$ follows from the last theorem.}
\label{proj3rmk}
\end{remark}

\begin{theorem} Let $M$ be an almost complex manifold of $\dim_{\C}\,M=2$.
Then $M$ has $(D_{c})$ ($\Rightarrow (D_{o})\Rightarrow (D_{r})$).
(This is in contrast with the results of Section 3 on surfaces in the algebraic
setting.)
\label{dim4}
\end{theorem}
\begin{proof}
We appeal to Lemma \ref{lemma1c}. The bundle $\F_{c}$ on $\boun X$ defined there
is a complex line bundle, and hence it extends to $X$ iff its first Chern class
$c_{1}(\F_{c})\in H^{2}(\boun X)$ lifts to $H^{2}(X)$. So it is enough to show that 
the restriction homomorphism
$$
H^{2}(X)\rarr H^{2}(\boun X)
$$
is surjective. From the commutative diagram induced by inclusions
\begin{eqnarray*}
\begin{array}{ccc}
H^{2}(X) &\rarr& H^{2}(\boun X)\\
j^{\ast}\uparrow& &\uparrow l^{\ast}\\
H^{2}(M\times M)&\rarr & H^{2}(U)
\end{array}
\end{eqnarray*}
we note that the bottom restriction map is the same as the map
$$
H^{2}(M\times M)\stackrel{\delta^{\ast}}{\rarr}H^{2}(M)
$$
(by deforming $U$ to its core $\Delta$) which is a split surjection. (Recall that $\delta$ denotes 
the diagonal embedding.)
The left vertical map is an isomorphism, for by excision
$$
H^{i}(M\times M,X) \cong H^{i}(D(\nu),S(\nu))
$$
and this vanishes for $0\leq i\leq 3$ by the Thom isomorphism ($\nu\cong \tau_{M}$ is a bundle
of real rank 4).
Also, the right vertical arrow is an isomorphism because of
$$
H^{i}(U,\boun X)=H^{i}(U,\boun U)= H^{i}(D(\nu),S(\nu))=0
$$
for $0\leq i\leq 3$, again by the Thom isomorphism.

\smallskip

Hence the top horizontal map is a surjection, and the assertion follows.
\end{proof}

\begin{lemma}[cf. \cite{BP}] Let $M$ be an almost complex manifold with $\dim_{\C}M=3$.
Let $\E$ be a smooth complex vector bundle on $M$ of any rank, with Chern classes
$c_{i}(\E)\in H^{2i}(M,\Z)$, for $1\leq i\leq 3$. Then these Chern classes satisfy
the following identity:
\begin{equation}
c_{3}(\E)-c_{2}(\E)(c_{1}(M)+c_{1}(\E))=2m\mu
\end{equation}
for some $m\in \Z$, where $\mu\in H^{6}(M,\Z)\cong \Z$ is a generator.
\label{cherneven}
\end{lemma}
\begin{proof}  See \cite{BP} for details. We remark that even though
this result is stated there for complex manifolds, the same proof works
by appealing to the generalised Riemann-Roch for almost complex
manifolds. Specifically, that the Todd characteristic class $T(M,\E)$ is integral (cf. 
the Theorem 25.5.4 in the Appendix 1 of \cite{Hirz}.)
\end{proof}

\bigskip

\begin{theorem} Let $M$ be an almost complex manifold of
$\dim_{\C}M=3$. Assume that $H^{1}(M,\Z)=0$ and $H^{2}(M,\Z)=\Z$. Then if
$M$ satisfies $(D_{c})$, the second Stiefel-Whitney class $w_{2}(M)$
vanishes (i.e. $M$ is spin).
\label{spinprop2}
\end{theorem}
\begin{proof}
Let $\E$ be a smooth complex rank 3 bundle on $M\times M$ realizing
$(D_{c})$. Let $x$ denote a generator of $H^{2}(M,\Z)$. Since $H^{1}(M,\Z)=0$, and
$H^{2}(M,\Z)=\Z x$, we have by
the K\"unneth formula that 
$$
H^{2}(M\times M)= \Z(x\times 1)\oplus \Z(1\times x)\,,
$$ 
where $\times$ denotes the cohomology cross product.
Hence the first  Chern class of $\E$ is given by
$$
c_{1}(\E) = a_{1}(x\times 1) + a_{2}(1\times x)\;\;\in H^{2}(M\times M,\Z)
$$
with $a_1, a_2\in \Z$. Since $\E$ restricted to the diagonal $\Delta$ is isomorphic as a {\it complex}
vector bundle to the normal bundle $\nu$ of $\Delta$, i.e. $\tau_{M}$, it follows that
$$
\delta^{\ast}(c_{1}(\E))=a_{1}(x.1)+a_{2}(1.x)=(a_{1}+a_{2})x=c_{1}(M)\,.
$$
Thus we have the relation (analogous to the ``weak point property'' from Section 4)
\begin{eqnarray}
(a_{1}+a_{2})x= c_{1}(M)\,.
\label{a1a2}
\end{eqnarray}

We noted in Remark \ref{euler1}, that the restriction of $\E$ to the slices $M\times
\{p\}$ and $\{q\}\times M$ will have Euler class $\pm 1$ times the
generator. Thus the 3rd Chern classes of the rank 3 bundles $\E_{1}:=\E_{|M\times \{p\}}$
and $\E_{2}:=\E_{|\{q\}\times M}$ are both equal to $\pm \mu$, where $\mu$ is the fundamental
class in $H^{6}(M,\Z)$. Clearly $c_{1}(\E_{1})=  a_{1}x$ and $c_{1}(\E_{2})=a_{2}x$.
From Lemma \ref{cherneven} applied to $\E_{1}$ and Eq. (\ref{a1a2}), it follows that
$$
c_{2}(\E_{1})(c_{1}(M)+c_{1}(\E_{1})) =  c_{2}(\E_{1})(2a_{1}+
a_{2})x= c_{3}(\E_{1}) + 2m_{1}\mu = (2m_{1}\pm 1)\mu
$$
and similarly
$$
c_{2}(\E_{2})(2a_{2} + a_{1})x= (2m_{2}\pm 1)
$$
for some $m_{i}\in \Z$. Reading these relations in $\Z_{2}$-cohomology, we find
$$
a_{1}\equiv a_{2}\equiv 1 \ \ \ \hbox{mod 2}\,.
$$
Thus $a_{1}+a_{2}\equiv 0$ mod 2. By Eq.~(\ref{a1a2}), it follows that modulo 2
$$
w_{2}(M)=c_{1}(M)=0\,.
$$

The theorem has been proved.
\end{proof}

\begin{remark}{\rm The condition $H^{2}(M,\Z)=\Z$ cannot be dropped
in Theorem \ref{spinprop2}. E.g., we know that $(D_{c})$ holds
for $\C\bP^2\times \C\bP^1$ (by Examples \ref{prodexample} and \ref{projgras})
which is not spin.} This condition seems to be the topological analogue of the condition $\Pic(X)=\Z$
in the algebraic theory. We also note that the converse to Theorem
\ref{spinprop2} would also be false. It is known that $S^{6}$ is an
almost complex manifold, which is clearly spin, but cannot satisfy
$(D_{c})$ since it does not satisfy $(D_{r})$ by the Theorem \ref{prop1}
above. 
\label{spinproprmk}
\end{remark}

\begin{corollary} Let $M\subset \C\bP^{N}$ be a smooth projective variety of
$\dim_{\C}M=3$. Assume that $M$ is a strict complete intersection
(alternatively, a set-theoretically complete
intersection with $H_{1}(M,\Z)=0$).
Then $M$ has $(D_{c})$ only if $M$ is spin.
\label{compinter}
\end{corollary}
\begin{proof} If $M$ is a smooth projective
set-theoretically complete intersection of complex dimension 3,
it is known (cf. Cor. 7.6 on p. 149 of \cite{Har1})
that for such an $M$ we have
$$
H^{1,0}(M)=H^{0,1}(M)=0\,.
$$ 
By Hodge decomposition, it follows that $H^{1}(X,\Z)=0$.
The same result quoted above shows that 
$$
H^{0,2}(M)=H^{2,0}(M)=0 \ \ \  \hbox{and} \ \ \ H^{1,1}(M)=\C\,.
$$ 
It follows again that $H^{2}(M,\C)=\C$. Since by hypothesis
$H_{1}(X,\Z)=0$, it follows that $H^{2}(X,\Z)$ has no torsion, equals $\Z$, and is
generated by the hyperplane class. Similarly, for a strict smooth complete intersection
$M$ with $\dim_{\C}M\geq 3$, it is known that $H^{2}(M,\Z)=\Z$ by the 
Grothendieck-Lefschetz theorem
(cf., e.g., \cite{Har1}). The result now follows from Theorem 
\ref{spinprop2}.
\end{proof}

\begin{corollary} Let $M$ be a smooth strict complete intersection of
$\dim_{\C}M=3$ in $\C\bP^{n}$, with $M=X_{1}\cap \cdots \cap X_{n-3}$ with
$X_{i}$ smooth hypersurfaces of degree $d_{i}$. Then $M$ has
$(D_{c})$ only if 
\begin{equation}\label{nd}
n+1-\sum_{i}d_{i}
\end{equation}
is even. In particular, a smooth hypersurface $M$ in $\C\bP^{4}$ has $(D_{c})$ only if it is 
of odd degree. Thus a smooth quadric $3$-fold in $\C\bP^{4}$ does not have $(D_{c})$.
\label{hyperP4}
\end{corollary}
\begin{proof}
It is known that the first Chern class of the tangent bundle of $M$
is $(n+1-\sum_{i}d_{i})$ times the hyperplane class
(indeed, the normal bundle of $M$ in $\C\bP^{n}$ is $\oplus_{i}\,{\cal O}(d_{i})$).
Thus $c_{1}(M)$ is an even multiple of the hyperplane
class iff the number (\ref{nd}) is even. Thus $M$ is spin iff this number is even, and the previous
corollary implies the result.
\end{proof}

When $n=4$, and $M$ is a hypersurface, the number (\ref{nd}) is even iff $d$ is odd.
In particular, a smooth quadric in $\C\bP^{4}$ does not have $(D_{c})$ (compare with
Proposition \ref{quadric1} from Section 4).

We remark that since the quadrics of complex dimension 1 and 2 are respectively
$\C\bP^{1}$ and $\C\bP^{1}\times \C\bP^{1}$, they both satisfy
$(D_{c})$. The quadric of dimension $3$ does not satisfy $(D_{c})$ by Corollary \ref{hyperP4}.
This last fact generalizes to all smooth projective quadric hypersurfaces of odd complex
dimension $\geq 3$.

\medskip

Up to the end of the proof of Theorem \ref{oddquad}, we shall now write $\bP^n$ for $\C\bP^n$.

\begin{proposition}
Let $Q_{2m-1}\subset \bP^{2m}$ denote the
smooth odd-dimensional quadric hypersurface $V(X_{0}^{2}+\cdots +X_{2m}^{2})$, and let $m\geq 2$. Then the
integral cohomology ring of $Q_{2m-1}$ is given by
$$
H^{\ast}(Q_{2m-1})=\Z[x,y]/\langle x^{m}-2y, y^{2}\rangle\,,
$$
where $x:=c_{1}({\cal O}_{Q_{2m-1}}(1))$ is the generator of
$H^{2}(Q_{2m-1})$, and $y$ is the generator of $H^{2m}(Q_{2m-1})$.
In particular,
\begin{eqnarray*}
H^{2k+1}(Q_{2m-1}) &=& 0\;\;\;\mbox{for all}\;\;k\\
H^{2k}(Q_{2m-1}) &=& \Z x^{k}\;\;\mbox{for all}\;\;0\leq k\leq m-1\\
 &=& \Z x^{k-m}y\;\;\mbox{for all}\;\;m\leq k\leq 2m-1\,.
\end{eqnarray*}
\label{quadricZ}
\end{proposition}
\begin{proof}
This result is well-known, but we sketch the proof for completeness. There is an inclusion
$j:\bP^{m-1}\hookrightarrow Q_{2m-1}$\,,
where $\bP^{m-1}$ is the linear subspace of $\bP^{2m}$ defined
by
$$
\{[x_{0}:x_{1}:\cdots :x_{2m}]\in \bP^{2m}:
x_{0}+\root\of{-1}x_{1}=\cdots =x_{2m-2}+\root\of{-1}x_{2m-1}=x_{2m}=0\}\,.
$$
Letting $i:Q_{2m-1}\hookrightarrow \bP^{2m}$ denote the natural inclusion
we have the composite homomorphisms
\begin{eqnarray*}
H^{r}(\bP^{2m})\stackrel{i^{\ast}}{\rarr}H^{r}(Q_{2m-
1})\stackrel{j^{\ast}}{\rarr}H^{r}(\bP^{m-1})\,,\\
H_{r}(\bP^{m-1})\stackrel{j_{\ast}}{\rarr}H_{r}(Q_{2m-1})
\stackrel{i_{\ast}}{\rarr} H_{r}(\bP^{2m})
\end{eqnarray*}
which are isomorphisms for $0\leq r\leq 2m-2$, since $i\circ j$ is a linear inclusion.
It follows that $H^{r}(Q_{2m-1})=H_{r}(Q_{2m-1})=0$ for
$r$ odd and $0\leq r\leq 2m-2$. Furthermore $H^{2k}(Q_{2m-1})=\Z x^{k}$
for $0\leq k\leq m-1$. By Poincar\'e duality on $Q_{2m-1}$, we have
$H^{i}(Q_{2m-1})=0$ for all odd $i$, $0\leq i\leq 4m-2$. Similarly,
$$
j_{\ast}:H_{2k}(\bP^{m-1})\rarr H_{2k}(Q_{2m-1})
$$
is an isomorphism of infinite cyclic groups for $0\leq k\leq m-1$.

\smallskip

Setting $D_{i}$, $i=1,2$, to be the Poincar\'e
duality isomorphisms for $\bP^{m-1}$ and $Q_{2m-1}$ respectively,
it follows by the preceding paragraph that the composition
\begin{eqnarray*}
H^{2k}(\bP^{m-1})\stackrel{D_{1}}{\rarr}H_{2m-2-2k}(\bP^{m-
1})\stackrel{j_{\ast}}{\rarr}H_{2m-2-2k}(Q_{2m-
1})\\
\stackrel{D_{2}^{-1}}
{\rarr} H^{2m+2k}(Q_{2m-1})
\end{eqnarray*}
is an isomorphism for $0\leq k\leq m-1$. This composite map is the
integral cohomology Gysin homomorphism denoted $j_{!}$, so, setting $k=0$, we find that
$H^{2m}(Q_{2m-1})$ is a cyclic group generated by $y=j_{!}1$. Also,
$j_{!}$ is a $H^{\ast}(Q_{2m-1})$-module homomorphism, so
$H^{2m+2k}(Q_{2m-1})$ is a cyclic group generated by
$$
j_{!}(h^{k})=x^{k}j_{!}1=x^{k}y
$$
for all $0\leq k\leq m-1$. Since $H^{4m}(Q_{2m-1})=0$, it follows that $y^{2}=0$. Since $Q_{2m-1}$ is a
degree 2 hypersurface in $\bP^{2m}$, we have
$$
\langle x^{2m-1},[Q]\rangle =2\,,$$
where $[Q]\in H_{4m-2}(Q_{2m-1})$ is the fundamental homology class of $Q_{2m-1}$. Thus
$$
\langle x^{m}.x^{m-1},[Q]\rangle =2\,.
$$
By Poincar\'e duality the generators $y$ of $H^{2m}$ and $x^{m-1}$ of $H^{2m-2}$ are dually paired,
so we have
$$
\langle y.x^{m-1},[Q]\rangle =1\,.
$$
Thus $x^{m}=2y$ \ and the proposition is proved.
\end{proof}

\begin{corollary}The cohomology ring $H^{\ast}(Q_{2m-1},\Z_{2})$ (where
$m\geq 2$) is given by
$$
H^{\ast}(Q_{2m-1},\Z_{2})=\Z_{2}[\xi,\eta]/\langle \xi^{m},
\eta^{2}\rangle\,,
$$
where $\xi$ (resp. $\eta$) is the mod 2 reduction of $x$ (resp. $y$) of the last
proposition. Alternatively, $\xi=w_{2}({\cal O}_{Q_{2m-1}}(1))$, the
second Stiefel-Whitney class of the canonical bundle on $Q_{2m-1}$ considered as a
real $2$-plane bundle, and $\eta=j_{!}1$, where
$$
j_{!}:H^{\ast}(\bP^{m-1},\Z_{2})\rarr H^{\ast+2m}(Q_{2m-1},\Z_{2})
$$
is the $\Z_{2}$-cohomology Gysin homomorphism. In particular,
\begin{eqnarray*}
H^{2k+1}(Q_{2m-1},\Z_{2})&=& 0\;\;\mbox{for all}\;\;k\,; \\
H^{2k}(Q_{2m-1},\Z_{2})&=& \Z_{2}\xi^{k}\;\;\mbox{for all}\;\;0\leq k\leq m-1\\
&=& \Z_{2}\xi^{k-m}\eta\;\;\mbox{for all}\;\;m\leq k\leq 2m-1\,.
\end{eqnarray*}
\label{Z2coho}
\end{corollary}
\begin{proof}
The assertion is immediate from the last proposition. For a complex vector bundle, the total
Stiefel-Whitney class is the mod 2 reduction of the total Chern class, so
$\xi=w_{2}({\cal O}_{Q_{2m-1}}(1))$.
\end{proof}

\begin{lemma}Let $m\geq 2$. Then the second Steenrod squaring operation $Sq^{2}$ on
$H^{\ast}(Q_{2m-1},\Z_{2})$ satisfies the relations
\begin{eqnarray*}
Sq^{2}(\xi) &=& \xi^{2}\\
Sq^{2}(\eta) &=& (m-1)\,\xi\eta \;\;\;\;\mbox{mod}\;\;2\\
Sq^{2}(\xi^{m-2}\eta)&=&\xi^{m-1}\eta\;\;\;\;\mbox{mod}\;\;2\,,
\end{eqnarray*}
where $\xi$ and $\eta$ are the algebra generators from Corollary \ref{Z2coho}.
\label{steenrod2}
\end{lemma}
\begin{proof}
Since $Sq^{i}x =x^{2}$ for $x\in H^{i}$ (cf. \cite{MS}, part (3) on  p.~90), and $\xi\in H^{2}(Q_{2m-1},\Z_{2})$,
it follows that $Sq^{2}(\xi) =\xi^{2}$.

For the second formula, one notes that the Gysin homomorphism $j_{!}$ is well known to be the composition
\begin{eqnarray*}
&&H^{i}(\bP^{m-1},\Z_{2})\stackrel{\phi}{\rarr}H^{i + 2m}(D(\nu),
S(\nu);\Z_{2})
\\
& &\stackrel{(l^{\ast})^{-1}}{\rarr}
 H^{i+2m}
(Q_{2m-1}, Q_{2m-1}\setminus {\mathbb P}^{m-1};\Z_{2})
\rarr H^{i+2m}(Q_{2m-1},\Z_{2})\,,
\end{eqnarray*}
where $\nu$ is the real rank $2m$ normal bundle of $\bP^{m-1}$ in $Q_{2m-1}$, $D(\nu)$ its disc bundle,
$S(\nu)$ its sphere bundle, $\phi$ the $\Z_{2}$ Thom isomorphism for $\nu$, $(l^{\ast})^{-1}$
is an excision isomorphism, and the last arrow is restriction. For brevity's sake, denote the composite of
the last two maps by $\alp$. Then
$$
\eta=j_{!}1=\alp(\phi(1))=\alp(U_{\nu})\,,$$
where $U_{\nu}\in H^{2m}(D(\nu), S(\nu);\Z_{2})$ is the $\Z_{2}$ Thom class of $\nu$.
Since $\alp$ is the composite of maps induced by restriction (and the inverse of a restriction),
the functorial operation $Sq^{2}$ commutes with $\alp$. Thus
$$
Sq^{2}(\eta)=Sq^{2}(\alp (U_{\nu}))=\alp(Sq^{2}U_{\nu})\,.
$$
It remains to determine $Sq^{2}U_{\nu}$. By Thom's identity for Stiefel-Whitney
classes (p.~91, {\em (loc.cit.)}), we have $\phi(w_{i}(\E))=Sq^{i}U_{\E}$ for any real
bundle $\E$, so $Sq^{2}U_{\nu}=\phi(w_{2}(\nu))$. The normal
bundle of $Q_{2m-1}$ in $\bP^{2m}$ is ${\cal O}_{Q_{2m-1}}(2)$,
and the normal bundle of the linear subspace
$\bP^{m-1}$ in $\bP^{2m}$ is the sum of
$m+1$ copies of ${\cal O}_{\bP^{m-1}}(1)$, hence
$$
\nu\oplus {\cal O}_{\bP^{m-1}}(2)=\left[{\cal O}_{\bP^{m-1}}(1)\right]^{m+1}\,.
$$
Thus
$$
c_{1}(\nu)=(m+1)h - 2h = (m-1)h\,,
$$
where $h$ is the hyperplane class of $\bP^{m-1}$. Since $w_{2}(\nu)$ is the mod 2 reduction of
$c_{1}(\nu)$, it follows that modulo 2
$$
Sq^{2}(U_{\nu})=\phi(w_{2}(\nu))=(m-1)\phi (h)\,.
$$
Thus
\begin{eqnarray*}
Sq^{2}(\eta)&=&\alp(Sq^{2}U_{\nu})=(m-1)\alp(\phi(h))\;\;\mbox{mod}\;\;2\\
&=&(m-1)j_{!}(h)=(m-1)\xi j_{!}(1)=(m-1)\xi \eta \;\;\mbox{mod}\;\;2
\end{eqnarray*}
since $j^{\ast}\xi=h$ mod 2 and $j_{!}$ is a $H^{\ast}(Q_{2m-1},\Z_{2})$-module homomorphism. This
proves the second formula.
Since $$
Sq^{k}(a.b)=\sum_{i+j=k}Sq^{i}(a)Sq^{j}(b)
$$
(cf. (4), p. 91, {\em (loc.cit.)}), and $Sq^{1}\equiv 0$ (odd cohomology vanish), we have by
the first two relations above

\begin{eqnarray*}
Sq^{2}(\xi^{m-2}\eta) &=& Sq^{2}(\xi^{m-2})\eta + \xi^{m-2}Sq^{2}(\eta)
=(2m-3)\xi^{m-1}\eta \;\;\;\mbox{mod}\;\;2\\
&=& \xi^{m-1}\eta
\end{eqnarray*}
which proves the third formula.

\smallskip

The lemma has been proved.
\end{proof}

\begin{proposition} Let $\E$ be a continuous complex vector bundle of any rank on
$Q_{2m-1}$ (where $m\geq 2$). In terms of the generators of $H^{2}$, $H^{4m-4}$ and
$H^{4m-2}$ determined in Proposition \ref{quadricZ}, define its Chern numbers $c_{j}\in \Z$
by
$$
c_{2m-1}(\E)=c_{2m-1}(x^{m-1}y);\;\;c_{2m-2}(\E)=c_{2m-2}(x^{m-2}y);\;\;c_{1}(\E)=c_{1}x
$$
Then we have
$$
c_{2m-1}=c_{2m-2}(c_{1}+1)\;\;\;\mbox{mod}\;\;2 \,\,.
$$
\label{cherniden2}
\end{proposition}
\begin{proof}
For a complex bundle $\E$,
the Stiefel-Whitney class $w_{2j}(\E)$ is the mod 2
reduction of $c_{j}(\E)$ (and odd Stiefel-Whitney classes vanish), so that
$$
w_{4m-2}(\E)=c_{2m-1}\xi^{m-1}\eta;\;w_{4m-4}(\E)=c_{2m-2}\xi^{m-
2}\eta;\;w_{2}(\E)=c_{1}\xi\;(\mbox{mod}\;2)\,,
$$
where $\xi$ and $\eta$ are the mod 2 reductions of $x$ and $y$
respectively, as in Corollary \ref{Z2coho}. We recall Wu's formula for
Stiefel-Whitney classes of a real bundle (cf. Problem 8-B on p. 94, {\it (loc.cit.)}):
$$
Sq^{2}w_{n}=w_{2}w_{n} + (2-n)w_{1}w_{n+1} +\fr{(2-n)(2-n-1)}{2}w_{0}w_{n+2}\,.
$$
For a complex vector bundle $\E$, $w_{1}(\E)=0$, so
applying the last formula for $n=4m-4$, we have
\begin{eqnarray*}
Sq^{2}(w_{4m-4}(\E)) &=& w_{2}(\E)w_{4m-4}(\E)+ w_{4m-2}(\E)\,.
\end{eqnarray*}

Substituting from the first paragraph and Corollary \ref{steenrod2}, we have
\begin{eqnarray*}
Sq^{2}(c_{2m-2}.\xi^{m-2}\eta) &=& (c_{1}\xi).(c_{2m-2}\xi^{m-2}\eta)+
c_{2m-1}\xi^{m-1}\eta\;\;\mbox{mod}\;\;2\\
c_{2m-2}(\xi^{m-1}\eta) &=& (c_{1}c_{2m-2}+c_{2m-1})\xi^{m-
1}\eta\;\;\mbox{mod 2}
\end{eqnarray*}
which implies the proposition.
\end{proof}

\begin{theorem} Let $m\geq 2$. Then a smooth quadric hypersurface $Q_{2m-1}\subset
\bP^{2m}$ does not have $(D_{c})$.
\label{oddquad}
\end{theorem}
\begin{proof}
The proof proceeds exactly as in the proof of Theorem \ref{spinprop2}. Since the normal bundle
of $Q_{2m-1}$ in $\bP^{2m}$ is ${\cal O}(2)$, the first Chern class of $Q_{2m-1}$ is
$$
c_{1}(\tau_{Q_{2m-1}})=(2m+1)x-2x=(2m-1)x\,.
$$
Thus
$$
w_{2}(\tau_{Q_{2m-1}})=\xi \in H^{2}(Q_{2m-1},\Z_{2})\,.
$$

Now let $\E$ be a complex vector bundle of complex rank
$2m-1$ on $Q_{2m-1}\times Q_{2m-1}$ realizing $(D_{c})$.
Then
$$
c_{1}(\E)=a_{1}(x\times 1) + a_{2}(1\times x)\in H^{2}(Q_{2m-1}\times Q_{2m-1},\Z)\,.
$$
Since $(D_{c})$ implies that $\delta^{\ast}(\E)\cong \tau_{Q_{2m-1}}$, we have
$$
\delta^{\ast}(c_{1}(\E))=
(a_{1}+a_{2})x=c_{1}(\tau_{Q_{2m-1}})=(2m-1)x\,,
$$
so that $a_{1}+a_{2}\equiv 1\;\;\mbox{mod}\;2$. The restrictions $\E_{i},\;\;i=1,2$ of $\E$
to the slices $Q_{2m-1}\times \{p\}$
and $\{q\}\times Q_{2m-1}$ respectively must have top Chern number
$c_{2m-1}\equiv 1\;\;\mbox{mod}\;2$ by Remark \ref{euler1}. This implies,
by Proposition \ref{cherniden2}, that $a_{1}$ and $a_{2}$ are both $\equiv 0$ mod 2. This
contradicts the last paragraph. The theorem now follows.
\end{proof}

\begin{remark}{\rm It is not clear what happens for quadrics of even complex dimension.
We note that $Q_{2}=\C\bP^{1}\times
\C\bP^{1}$ and $Q_{4}=G_{2}(\C^{4})$ both satisfy $(D_{c})$ by Examples \ref{prodexample} and
\ref{projgras}.}
\end{remark}

We now make use of the topological Theorem \ref{oddquad}, or rather its
key input, Proposition~\ref{cherniden2}, to obtain the following algebraic
result.
\begin{theorem}\label{oddquad-alg}
 Let $X\subset\bP^{2n}$ be an odd dimensional smooth
quadric hypersurface over an algebraically closed field $k$, with $n>1$.
Then $X$ does not have (D).
\end{theorem}
\begin{proof}
We claim that the congruence formula in Proposition~\ref{cherniden2}
is valid for the algebraic Chern classes of any algebraic vector bundle,
if we consider it as taking values in the Chow ring $CH^*(X)\tensor\Z/2\Z$.
Assuming the claim, it follows that $X$ does not have the $\O_X(-r)$-point
property for any odd $r>0$, exactly as in the complex topological case.
Since $\omega_X=\O_X(-(2n-1))$, we see from Corollary~\ref{Lpt} that $X$
does not have (D).

Indeed, we will show that the congruence actually holds for Chern classes
of virtual bundles, i.e. for the homogeneous components of any element in
the image of the mod 2 total Chern class map
$$
K_0(X)\to CH^*(X)\tensor\Z/2\Z\,.
$$

We first observe that any smooth quadric $X\subset \bP^{2n}$ over an
algebraically closed field $k$ is isomorphic to the quadric defined by
$$
\{x_0^2+\sum_{i=1}^nx_ix_{i+n}=0\}.
$$ 
This quadric contains the
$i$-dimensional projective linear subspace
$$L_i,\;\; 0\leq i\leq n-1,$$
defined by
$$x_0=0,x_j=x_{n+j}=0\;\forall\;1\leq
j<n-i,\;\;x_j=0\,\forall\,j\geq n-i.$$

Here $L_{n-1}\cong \bP^{n-1}$ is defined by $x_0=x_1=\ldots=x_n=0$, and
$L_0\subset L_1\subset\cdots\subset L_{n-1}$ is a maximal flag of
projective linear subspaces. Let $H_j\subset X$ be a complete intersection
with a general linear subspace of $\bP^{2n}$ codimension $j$.

The Grothendieck group $K_0(X)$ of algebraic vector bundles on $X$ is
freely generated by the classes $\O_{H_j}$, $0\leq j\leq n-1$, and
$\O_{L_i}$, $0\leq i\leq n-1$. Here $H_0=X$.

Similarly the Chow ring of $X$ is freely generated, as an abelian group,
by the classes of the irreducible varieties $H_j$, $0\leq j\leq n-1$ and
$L_i$, $0\leq i\leq n-1$. Here $[H_j]\in CH^j(X)$, while $[L_j]\in
CH_j(X)=CH^{2n-1-j}(X)$.

If $x\in CH^1(X)$ is the class of the cycle $[H_1]$, and $y\in CH^n(X)$
the class  of the cycle $[L_{n-1}]$, then the intersection product in
$CH^*(X)$ is determined by the properties that the class of $[H_j]$ is
just $x^j$. and the class of $[L_j]$ is $x^{n-j-1}y$. We also have
relations $y^2=0$, and $x^n=2y$.

For $k=\C$, this is exactly the presentation for the integral cohomology
ring $H^*(X,\Z)$ given in Proposition~\ref{quadricZ} (with $m$ in place
of $n$). In fact, the cycle class map
\[CH^*(X)\to H^*(X,\Z)\]
induces an isomorphism of rings, where the algebraic generators $x,y$
above map to the corresponding topological generators of
Proposition~\ref{quadricZ}.

In a similar fashion, if $K_{top}(X)$ denotes the topological $K$-group of
complex vector bundles on $X$,  one knows that the natural map $K_0(X)\to
K_{top}(X)$ is an isomorphism of rings. Further, the algebraic and
topological Chern class maps are compatible. Hence, an element of
$CH^*(X)$ lies in the image of the algebraic total Chern class map
precisely when the corresponding element in $H^*(X,\Z)$ is in the image of
the topological total Chern class map. Note that, at the topological
level, any element of $K_{top}(X)$ is of the form $[\E]-m$ for some integer
$m$, where $\E$ is a complex topological vector bundle on $X$. Thus the
image of the total topological Chern class map coincides with the set of
total Chern classes of complex topological vector bundles. In particular,
at the topological level, the congruence in Proposition~\ref{cherniden2}
is valid for the total Chern class of any element on $K_{top}(X)$. Hence
this congruence is valid for the corresponding algebraic total Chern class
image of $K_0(X)$ in $CH^*(X)$.

We now argue that, in fact, the image of the algebraic total Chern
class map is the same, for any smooth quadric $X\subset\bP^{2n}$ over
an algebraically closed field $k$, where we identify the Chow rings of
all such quadrics using the explicit presentation. This is because the
formulas for the Chern classes of the generators $[\O_{H_j}]$ and
$[\O_{L_i}]$ of  $K_0(X)$, expressed as polynomials over $\Z$ in the
generators $x,y$ for the Chow ring, are independent of the field $k$.
This is clear for the $H_j$ since these are classes pulled back from
projective space, and Chern classes are functorial under pullback; these
classes are certain polynomials in $x$. For the classes of $L_i$, we may
regard these as in the image of the  push-forward map $K_0(L_{n-1})\to
K_0(X)$. The push-forward map on the graded Chow group is also defined,
identifying $CH^j(L_{n-1})\cong \Z$ with  $CH^{n+j}(X)\cong \Z$, where
the cycle $[L_j]$ is a generator of either group.

Now the universal nature of the Chern class formulas for the
sheaves $\O_{L_i}$ follows from the theorem of Riemann-Roch ``without
denominators'' (\cite{F2}, Theorem~15.3) for this push-forward map.

Thus, the truth of the congruences of Proposition~\ref{cherniden2} in the
topological case, hence the complex algebraic case, implies the same
congruences hold over any algebraically closed field.
\end{proof}
\begin{remark} In fact the Riemann-Roch formula in \cite{F2}, Theorem~15.3
allows us, in principle, to consider giving a direct algebraic proof of
Proposition~\ref{cherniden2} for elements in the image of the algebraic
total Chern class map, by actually determining explicitly the image of
the total Chern class map itself. This would also prove the result
in the topological case as well, by reversing the above argument,
and would in some way ``explain'' it, independent of any identities or
properties of Steenrod squares, for example. However, we were so far
unable to make the direct computation. In principle, such a direct
approach might also shed light on the case of even dimensional quadrics,
and perhaps other cases. We conjecture that a smooth quadric hypersurface
$Q_n$ in the projective $(n+1)$-space (over an algebraically closed field)
has (D) iff $n=1,2$ or $4$.
\end{remark}

\smallskip

\noindent
{\bf Note} \ In a recent letter \cite{Deb}, O. Debarre has informed us that he
proved (D) for the Jacobians of curves, and found counterexamples to (D) for
general abelian varieties of dimension greater than 2.

\bigskip

{\bf Acknowledgments} \ The first two authors started cooperation on this
paper during the Conference ``Antalya Algebra Days 2005'', and they thank
the organizers for creating this opportunity.

P. Pragacz thanks Professor Friedrich Hirzebruch for his interest towards
the ``diagonal property'' as well as for helpful comments received before
and during the author's stay at the MPIM in Bonn in 2006.

V. Pati is grateful to the other two authors for inviting him on
board this interesting investigation, and valuable insights. Also to
Jishnu Biswas for educating him on complete intersections and pointing out
the relevant references, and Shreedhar Inamdar for useful discussions and
encouragement.

The authors also thank A. Ozgur Kisisel for pointing out an error in a former
version of the paper.

\end{document}